\documentclass{article}

\usepackage[margin=1in]{geometry}

\usepackage{graphicx}

\usepackage[utf8]{inputenc} 
\usepackage[T1]{fontenc}    
\usepackage{hyperref}       
\usepackage{url}            
\usepackage{booktabs}       
\usepackage{amsfonts}       
\usepackage{nicefrac}       
\usepackage{microtype}      

\usepackage{bm}
\usepackage{amsmath,amsthm}
\usepackage{mathtools,nccmath}
\usepackage{algorithmic}

\usepackage[algoruled,boxed,lined]{algorithm2e}
\usepackage{xcolor}         
\usepackage{epstopdf}

\usepackage{multirow}

\usepackage{enumitem}

\newcommand{\bx}{\bm{x}}

\newcommand{\by}{\bm{y}}

\newcommand{\bs}{\bm{s}}
\newcommand{\bg}{\bm{g}}
\newcommand{\bv}{\bm{v}}

\newcommand{\bu}{\bm{u}}
\newcommand{\bh}{\bm{h}}

\newcommand{\R}{\mathbb{R}}

\newcommand{\sO}{\mathcal{O}}

\newcommand{\G}{\mathsf{G}}

\usepackage{subcaption}
\usepackage{verbatim}
\newtheorem{assumption}{Assumption}
\newtheorem{condition}{Condition}

\newtheorem{lemma}{Lemma}

\theoremstyle{remark}

\newtheorem{remark}{Remark}

\usepackage{thmtools}
\usepackage{thm-restate}
\usepackage{float}

\makeatletter
\newcommand{\algline}[1]{%
  \begingroup
  \edef\@currentlabel{\theALC@line}%
  \label{#1}%
  \endgroup
}
\makeatother

\allowdisplaybreaks

\makeatother

\newcounter{relctr}[section] 

\AtBeginEnvironment{align}{\setcounter{relctr}{0}} 

\newcommand\labelrel[2]{%
  \begingroup
    \refstepcounter{relctr}%
    \stackrel{\textnormal{(\roman{relctr})}}{\mathstrut{#1}}%
    \originallabel{#2}%
  \endgroup
}

\AtBeginDocument{\let\originallabel\label}

\title{\textbf{Accelerating Single-Point Zeroth-Order Optimization with Regression-Based Gradient Surrogates}}

\author{%
   Xin Chen$^\dagger$\thanks{Corresponding author.} \\  Department of Electrical and Computer Engineering\\ 
   Texas A\&M University\\ 
   \texttt{xin\_chen@tamu.edu}\\ \\
   Zhaolin Ren\thanks{X. Chen and Z. Ren contributed equally to this work.} \\ Mitsubishi Electric Research Laboratories (MERL)\\ 
   \texttt{zren@merl.com} 
}

\date{}

\begin{document}

\maketitle
\begin{abstract}
 \fontsize{10pt}{13pt}\selectfont 
 Zeroth-order optimization (ZO) is widely used for solving black-box optimization and control problems. In particular, single-point ZO (SZO) is well-suited to online or dynamic problem settings due to its requirement of only a single function evaluation per iteration. However, SZO suffers from high gradient estimation variance and slow convergence, which severely limit its practical applicability. To overcome this limitation, we propose a novel yet simple SZO framework termed regression-based SZO (ReSZO), which substantially enhances the convergence rate. Specifically, ReSZO constructs a surrogate function via regression using historical function evaluations and employs the gradient of this surrogate function for iterative updates. Two instantiations of ReSZO, which fit linear and quadratic surrogate functions respectively, are introduced. Moreover, we provide a non-asymptotic convergence analysis for the linear instantiation of ReSZO, showing that its convergence rates are comparable to those of two-point ZO methods. Extensive numerical experiments demonstrate that ReSZO empirically converges two to three times faster than two-point ZO in terms of function query complexity.
\end{abstract}

\fontsize{11pt}{14pt}\selectfont

\section{Introduction}

In this paper, we consider solving a generic black-box unconstrained optimization problem \eqref{eq:problem}:
\begin{align}\label{eq:problem}
    \min_{\bx\in \R^d} \, f( \bx),
\end{align}
where $\bx\in \R^d$ is the decision variable of dimension $d$. The objective function $f:\R^d\to \R$ is unknown, and only its zeroth-order information, namely function evaluations, is accessible. As a result, the first-order gradients of function 
$f$ are unavailable, and problem \eqref{eq:problem} can be solved solely based on function evaluations. 
This black-box optimization has attracted considerable recent attention and covers a broad range of applications, including power system control \cite{chen2021safe,chen2020model}, reinforcement learning \cite{malik2019derivative,li2019distributed},  neural network training \cite{chen2023deepzero,chen2017zoo}, sensor management \cite{liu2020primer}, etc.

Zeroth-order optimization (ZO) methods \cite{liu2020primer,nesterov2017random} have been extensively developed and applied to solve black-box optimization problems, offering theoretical convergence guarantees. Essentially, ZO methods estimate unknown gradients using perturbed function evaluations. According to the number of queried function evaluations per iteration, ZO methods can be categorized into single-point, two-point, and multi-point ZO schemes \cite{liu2020primer}. Single-point ZO (SZO) methods \cite{flaxman2005online,zhang2021new,chen2022improve} require only a single function evaluation per iteration for gradient estimation, making them well-suited for online optimization or dynamic control applications. Because in such settings, the system environment is non-stationary, or each function evaluation involves implementing a decision that directly changes the system environment, rendering multiple function queries at the same instant infeasible. However, single-point ZO suffers from large estimation variance and slow convergence, severely limiting its practical applicability. Despite recent advancements \cite{zhang2021new,chen2022improve,saha2011improved,dekel2015bandit,gasnikov2017stochastic}, the convergence of single-point ZO remains slower than that of two-point ZO or multi-point ZO, which use multiple function evaluations for gradient estimation per iteration and achieve lower variance.

\textbf{Motivation}. At each iteration $t$,
classic ZO methods estimate the gradient $\nabla f(\bx_t)$ solely based on the function evaluations at the current iterate $\bx_t$, while all historical function evaluations obtained from previous iterations are discarded.  
However, these past function evaluations contain valuable information about the unknown function and can help with the gradient estimations for subsequent iterations, especially those obtained in the most recent iterations.
Moreover, in many practical applications, e.g., clinical trials, real-world experiments, and long-horizon simulations, obtaining even a single function evaluation can be costly and laborious, making it crucial to fully exploit all available evaluations to improve efficiency. Then, a natural question is 

\vspace{3pt}
``\emph{Can available historical function evaluations from previous iterations be leveraged to enhance single-point zeroth-order gradient estimation?}"
\vspace{3pt}

\textbf{Contribution}. Motivated by this question, we propose a novel yet simple single-point ZO framework termed  \textbf{Re}gression-based \textbf{S}ingle-point \textbf{Z}eroth-order \textbf{O}ptimization (\textbf{ReSZO}). Unlike classical ZO methods that estimate gradients using only function evaluations at the current iterate, ReSZO exploits recent historical function evaluations to significantly improve gradient estimation  while requiring only one function query per iteration. Our main contributions are summarized as follows: \vspace{-5pt}
\begin{itemize}
    \item [1)] (\emph{New Framework and Algorithms}). We introduce the ReSZO framework and develop two concrete instantiations, namely L-ReSZO and Q-ReSZO, which construct linear and quadratic surrogate models, respectively. By leveraging historical function evaluations together with a new adaptive perturbation scheme, the proposed methods substantially improve convergence efficiency over existing SZO methods and demonstrate strong empirical performance.
   \item [2)] (\emph{Theoretical Analysis}). We provide  theoretical analysis for the convergence of ReSZO. Specifically, for both smooth nonconvex and strongly convex objectives, we show that under certain conditions, 
   L-ReSZO achieves convergence rates that are within a factor \(C_d \!>\! 0\) of those attained by two-point ZO methods (see Table~\ref{tab:zo_complexity}). While tightly bounding \(C_d\) is challenging, we offer theoretical justification that $C_d$ scales as $\mathcal{O}(d)$. We also provide extensive empirical evidence indicating that, under suitable algorithmic settings, the practical scaling follows \(C_d = \mathcal{O}(\sqrt{d})\) (see Section~\ref{sec:analysis}).


   
    
    \item  [3)] (\emph{Empirical Evaluation}). We conduct extensive numerical experiments on a variety of benchmark problems, demonstrating that ReSZO not only matches but empirically
converges two to three times faster than two-point ZO in terms of function query complexity.  
\end{itemize}

\begin{table*}
  \captionsetup{width=\textwidth}
  \centering
  \caption{Iteration complexity of ZO methods for solving Lipschitz
           and smooth objective functions.}
  \label{tab:zo_complexity}

  \begin{tabular*}{\textwidth}{@{\extracolsep{\fill}} l l c c @{}}
    \toprule[1.2pt]
    \multirow{2}{*}{{Method}} &
    \multirow{2}{*}{{Literature}} &
    \multicolumn{2}{c}{{Complexity}} \\ \cmidrule(r){3-4}
    & & $\mu$-{Strongly convex} & {Non-convex} \\
    \midrule
   Classic SZO & \cite{gasnikov2017stochastic} &
      $\mathcal{O}\!\bigl(d^{2}/\varepsilon^{3}\bigr)$ & -- \\ 
      Residual-feedback SZO & \cite{zhang2022new} &
      $\mathcal{O}\!\bigl(d^{3}/\varepsilon^{\frac{3}{2}}\bigr)$ &
      $\mathcal{O}\!\bigl(d^{3}/\varepsilon^{\frac{3}{2}}\bigr)$ \\
    HLF-SZO & \cite{chen2022improve} &
      $\mathcal{O}\!\bigl(d^{\frac{3}{2}}/\varepsilon^{\frac{3}{2}}\bigr)$ &
      $\mathcal{O}\!\bigl(d^{\frac{3}{2}}/\varepsilon^{\frac{3}{2}}\bigr)$ \\ 
      \textbf{L-ReSZO} & \textbf{\textsc{This work}} &
      $\mathcal{O}\!\Bigl(C_d\frac{d}{\mu}\log\!\frac{1}{\varepsilon}\Bigr)$ &
      $\mathcal{O}\!\bigl(C_dd/\varepsilon\bigr)$ \\
    Two-point ZO & \cite{nesterov2017random} &
      $\mathcal{O}\!\Bigl(\frac{d}{\mu}\log\!\frac{1}{\varepsilon}\Bigr)$ &
      $\mathcal{O}\!\bigl(d/\varepsilon\bigr)$ \\ 
    \bottomrule[1.2pt]
  \end{tabular*}

  \vspace{-0.1em}
  \begin{flushleft}
\footnotesize\itshape
{\normalfont *In the strongly convex case, the accuracy is measured by
$\,f(x_{T})-f^*\le\varepsilon$, where
$f^* := \min_{x\in\mathbb{R}^{d}}f(x)$. 
In the non-convex case, the accuracy is measured by
$\tfrac1T\sum_{k=1}^{T}\lVert\nabla f(x_{k})\rVert^{2}\le\varepsilon$. 
}
\end{flushleft}
\end{table*}

\textbf{Related Work}. The classic SZO method, formulated in \eqref{eq:szo}, is first proposed by  \cite{flaxman2005online}, which estimates gradients by querying only a single function value per iteration.
In \cite{dekel2015bandit}, a self-concordant barrier function is used to enhance the convergence of SZO for solving bandit smooth convex optimization. 
By leveraging high-order smoothness, reference \cite{novitskii2021improved} derives an improved complexity bound for SZO in the context of stochastic convex optimization. 
To further enhance the convergence of SZO, reference \cite{zhang2022new} proposes a residual-feedback SZO method, formulated in \eqref{eq:rfszo}, which reuses the function evaluation from the last iteration for gradient estimation, resulting in significant performance improvement. Additionally, reference \cite{chen2022improve} leverages the idea of  high-pass and low-pass filters from extremum seeking control to accelerate SZO convergence. Notably, the integration of a high-pass filter coincides with the residual-feedback approach in \cite{zhang2022new}. This paper follows a similar direction to \cite{chen2022improve,zhang2022new}; however, rather than reusing only the most recent function evaluation, we leverage many historical function evaluations from previous iterations to fit surrogate models. In addition, our ReSZO method is related to the linear regression and interpolation approaches \cite{conn2009introduction,berahas2022theoretical,wang2024relizo} used in derive-free optimization (DFO) \cite{larson2019derivative}. The fundamental difference is that these methods \cite{conn2009introduction,berahas2022theoretical,wang2024relizo} perform multiple function evaluations near the current iterate for linear interpolation or regression. In contrast, our method queries the function value only once at each iteration and fits surrogate models by  recycling historical function evaluations. A recent SZO algorithm~\cite{jongeneel2024small} achieves strong convergence rates; however, unlike our approach, it assumes access to an oracle that supports complex-valued queries, which limits its applicability in standard zeroth-order settings where only real-valued function evaluations are available.




\section{Preliminaries} 

The classic single-point ZO (SZO) method \cite{flaxman2005online} that solves \eqref{eq:problem} is given by \eqref{eq:szo}:
\begin{align}\label{eq:szo}
  \ \textbf{Single-point ZO (SZO)}: 
  \ \bx_{t+1} = \bx_t - \eta\cdot \underbrace{\frac{d}{\delta} f(\bx_t+ \delta \bu_t) \bu_t}_{:=\, \G^{1}_f(\bx_t;\bu_t,\delta)},\quad t=0,1,2,\cdots
\end{align}
where $\delta\!>\!0$ is a scalar parameter called smoothing radius, $\eta \!>\!0$ is the step size, and $\bu_t\in\R^d$ is an independent identically distributed random direction vector sampled from $\mathrm{Unif}(\mathbb{S}_{d-1})$, namely the uniform distribution on the unit sphere $\mathbb{S}_{d-1}$, for all $t=0,1,2,\cdots$. The random direction vector $\bu_t$ can also be sampled from the multivariate Gaussian distribution $\mathcal{N}(\mathbf{0},\bm{I})$ or other zero-mean distributions. In \eqref{eq:szo}, $\G^{1}_f(\bx_t;\bu_t,\delta)$ denotes the SZO gradient estimator of function $f$.

Recent works have introduced a new variant of SZO, referred to as residual-feedback SZO (RSZO) \cite{zhang2022new} or high-pass filter-integrated SZO \cite{chen2022improve}, which is given by \eqref{eq:rfszo}:
\begin{align}\label{eq:rfszo}
 & \ \textbf{Residual-feedback SZO (RSZO)}: \text{for } t=1,2,\cdots,  \nonumber \\    
 & \ \bx_{t+1} = \bx_t \!- \!\eta\cdot \underbrace{\frac{d}{\delta}\Big( f(\bx_t\!+\! \delta \bu_t)\!-\!f(\bx_{t-1}\!+\! \delta \bu_{t-1})\Big) \bu_t}_{ \coloneq\, \G^{\mathrm{r}1}_f(\bx_t;\bu_t,\delta)}.
\end{align}
At each iteration $t$, the gradient estimator $\G^{\mathrm{r}1}_f(\bx_t;\bu_t,\delta)$ of
the RSZO method \eqref{eq:rfszo} reuses the function evaluation $f(\bx_{t\!-\!1}\!+\! r \bu_{t\!-\!1})$ queried at the last iteration $t\!-\!1$. This approach substantially  reduces gradient estimation variance and improves convergence compared with SZO \eqref{eq:szo}.

The classic two-point ZO (TZO) method \cite{nesterov2017random} that solves problem \eqref{eq:problem} is given by \eqref{eq:tzo}:
\begin{align}
    & \  \textbf{Two-point ZO (TZO)}: \text{ for } t=0,1,2,\cdots  \nonumber\\
   & \ \bx_{t+1} = \bx_t  - \eta\cdot \underbrace{\frac{d}{2\delta} \Big(f(\bx_t+ \delta \bu_t) -f(\bx_t-\delta\bu_t)\Big)\bu_t}_{ \coloneq\, \G^{2}_f(\bx_t;\bu_t,\delta)},
   \label{eq:tzo}
\end{align}
which employs two function evaluations to estimate the gradient of function $f$ with the two-point gradient estimator $\G^{2}_f(\bx_t;\bu_t,\delta)$.

The classic ZO gradient estimators above share the same expectation \cite{nesterov2017random}:
\begin{align}
    &  \mathbb{E}_{\bu\sim \mathrm{Unif}(\mathbb{S}_{d-1})}\big[\G^{1}_f(\bx;\bu,\delta)\big] \nonumber\\
    = \, & \mathbb{E}_{\bu\sim \mathrm{Unif}(\mathbb{S}_{d-1})}\big[\G^{\mathrm{r}1}_f(\bx;\bu,\delta)\big]\nonumber \\
  = \, &  \mathbb{E}_{\bu\sim \mathrm{Unif}(\mathbb{S}_{d-1})}\big[\G^{2}_f(\bx;\bu,\delta)\big] = \nabla \hat{f}(\bx;\delta), \label{eq:expect}
\end{align}
where $\hat{f}(\bx;\delta) := \mathbb{E}_{\bv\sim \mathrm{Unif}(\mathbb{B}_d)}\big[ f(\bx+\delta \bv)\big]$ is a smoothed version of function $f$ \cite{flaxman2005online}. 
Despite having the same expectation, the two-point ZO gradient estimator $\G^{2}_f(\bx_t;\bu_t,\delta)$ exhibits lower variance and achieves faster convergence, compared with 
the single-point ZO gradient estimators $\G^{1}_f(\bx_t;\bu_t,\delta)$ and $\G^{\mathrm{r}1}_f(\bx_t;\bu_t,\delta)$.
 See \cite{nesterov2017random,chen2022improve} for more details. 
 

\section{Regression-Based Single-Point Zeroth-Order Optimization} \label{sec:method}


In this paper, we focus on the class of SZO methods that require only one function evaluation per iteration. Below, we introduce the ReSZO method and its two instantiations, L-ReSZO and Q-ReSZO, which fit linear and quadratic surrogate models using historical function evaluations. 

 At each iteration $t$, the available function evaluations from all previous iterations are given by $\{f(\hat{\bx}_k)\}_{k=0}^{t}$ with $\hat{\bx}_k = \bx_k+\delta\bu_k$. Here, $\bu_k\in\R^d$ is a random direction vector sampled
from the uniform distribution $\mathrm{Unif}(\mathbb{S}_{d-1})$ over the unit sphere $\mathbb{S}_{d-1}$, 
and $\delta>0$ is a scalar parameter called smoothing radius. These are the same as the classic ZO methods \eqref{eq:szo}-\eqref{eq:tzo}.
Note that function evaluations are performed at the perturbed points $\hat{\bx}_k$ rather than the iterate points ${\bx}_k$ to promote exploration of the optimization space. The necessity of introducing such perturbations is further illustrated by 
Remark \ref{remark:pert} and the numerical experiments in Section \ref{sec:pert}.

\subsection{Linear Regression-Based Single-Point Zeroth-Order Optimization (L-ReSZO)}

Based on the first-order Taylor expansion at the current perturbed point $\hat{\bx}_t$, the available function evaluations $\{f(\hat{\bx}_k)\}_{k=0}^{t}$ can be expressed as \eqref{eq:1taylor}: for $k = 0,1\dots,t$,
\begin{align} \label{eq:1taylor}
f(\hat{\bx}_k)\! =\! f(\hat{\bx}_t) \!+\! \nabla f(\hat{\bx}_t)^\top \big(\hat{\bx}_k \!-\! \hat{\bx}_t\big)\! +\! \epsilon_{t,k}(\hat{\bx}_k \!-\! \hat{\bx}_t),
\end{align}
where the function $\epsilon_{t,k}(\cdot)$ represents the corresponding high-order residual term. Accordingly, we propose to fit a linear surrogate function $f^\mathrm{s1}_t(\bx)$ in \eqref{eq:sur1} to locally approximate the original function $f(\bx)$ in \eqref{eq:problem}: 
\begin{align}\label{eq:sur1}
    f^\mathrm{s1}_t(\bx) \coloneq   f(\hat{\bx}_t) +  \bg_t^\top \big({\bx} - \hat{\bx}_t\big),   
\end{align}
where the linear coefficient vector $\bg_t\in \R^d$ can be obtained through linear regression on the available sample points $\{\hat{\bm{x}}_{k}, f(\hat{\bm{x}}_{k})\}_{k=0}^t$. Then, the gradient of the surrogate function $\nabla  f^\mathrm{s1}_t(\bx) = \bm{g}_t$ can be used as the estimator for the true gradient $\nabla f(\bx)$. A key feature of the linear surrogate function \eqref{eq:sur1} is the omission of a constant intercept term, which enforces $f^\mathrm{s1}_t(\hat{\bx}_t) = f(\hat{\bx}_t)$ and reduces model bias.


To ensure that the fitted surrogate model $f^\mathrm{s1}_t(\bx)$ remains up-to-date with the most recent function evaluations, we employ recursive linear regression in a sliding window, which uses only the latest $m$ function evaluations $\{\hat{\bm{x}}_{t-i+1}, f(\hat{\bm{x}}_{t-i+1})\}_{i=1}^{m}$  to fit the linear model at each iteration $t$.
Specifically, at each iteration $t$, we solve the following least squares problem \eqref{eq:lr} to obtain $\bm{g}_t$:
\begin{align} \label{eq:lr}
     \ \min_{\bm{g}_t}\, \sum_{i=1}^{m}  \Big[  f^\mathrm{s1}_t(\hat{\bx}_{t-i+1})- f(\hat{\bx}_{t-i+1})\Big]^2  \Longrightarrow   \ 
   \bm{g}_t = (X_t^\top X_t)^{\dagger}X_t^\top  \bm{y}_t,
\end{align}
where $\dagger$ denotes the Moore–Penrose pseudoinverse, $\bm{\Delta}^x_{t,i}\!\coloneq\! \hat{\bx}_{t-i} \!-\! \hat{\bx}_t, \Delta^f_{t,i}\!\coloneq\! f(\hat{\bx}_{t-i}) \!-\! f(\hat{\bx}_{t})$, 
\begin{align}
X_t \coloneq \!\!\begin{bmatrix}
        {\bm{\Delta}^x_{t,1}}^\top \\  {\bm{\Delta}^x_{t,2}}^\top \\ \vdots \\  {\bm{\Delta}^x_{t,m-1}}^\top 
    \end{bmatrix}\!\!\in\R^{(m\!-\!1)\times d}, \ \bm{y}_t  \coloneq  \!\!\begin{bmatrix}
       \Delta^f_{t,1}\\\Delta^f_{t,2}\\ \vdots\\ \Delta^f_{t,m-1}
    \end{bmatrix}\!\!\in\R^{m\!-\!1}.
\end{align}
Consequently, the proposed \textbf{L}inear \textbf{Re}gression-based \textbf{S}ingle-point \textbf{Z}eroth-order \textbf{O}ptimization (\textbf{L-ReSZO}) method for solving problem \eqref{eq:problem} is presented as Algorithm \ref{alg:lReSZO}. Since the method requires $m$ function evaluations for regression, we initially employ a classical ZO algorithm, e.g., the RSZO method \eqref{eq:rfszo}, to perform iterative updates and accumulate sufficient sample points, as shown in Steps \ref{st:rszo:for}-\ref{st:rszo:end}. After $m$ iterations, the algorithm switches to using the linear regression-based gradient estimator $\bm{g}_t$ for subsequent iterative updates, as shown in Steps \ref{st:reszo:for}-\ref{st:reszo:end}. 
The sample size $m$ is typically chosen to be comparable to the dimension $d$ of $\bx$, or slightly larger, to ensure that  $X_t$ is of full column rank and the matrix $X_t^\top X_t$ in \eqref{eq:lr} is invertible. In this case, one can replace the pseudoinverse $(\cdot)^\dagger$ in \eqref{eq:lr} by the regular matrix inverse $(\cdot)^{-1}$.



\begin{algorithm}
\caption{Linear Regression-based Single-point Zeroth-order Optimization (L-ReSZO)}
\label{alg:lReSZO}
\begin{algorithmic}[1]
\STATE \textbf{Input}: Initial point $\bx_{0}$; step sizes $\tilde{\eta}, \eta$; smoothing radius $\tilde{\delta}$; sample size $m$; total iterations $T$. \algline{st:input}
\FOR{$t = 0$ to $m-1$} \algline{st:rszo:for}
    \STATE Perform RSZO update (i.e., Eq.~\eqref{eq:rfszo}) with $\tilde{\eta}, \tilde{\delta}$. \algline{st:rszo:update}
\ENDFOR \algline{st:rszo:end}

\FOR{$t = m$ to $T-1$} \algline{st:reszo:for}
    \STATE Use an \emph{adaptive smoothing radius} by setting: \algline{st:reszo:delta}
    \begin{align} \label{eq:adaptive}
        \delta_t \!=\! \| \bm{x}_{t}\!-\!\bm{x}_{t-1}\|.
    \end{align} 
    \STATE Sample $\bu_t \sim \mathrm{Unif}(\mathbb{S}_{d-1})$; set $\hat{\bx}_t = \bx_t + \delta_t \bu_t$;  evaluate $f(\hat{\bx}_t)$; and collect the latest $m$ samples $\{(\hat{\bx}_{t-i+1}, f(\hat{\bx}_{t-i+1}))\}_{i=1}^{m}$. \algline{st:reszo:func}
    \STATE Solve for the linear coefficient $\bg_t$ via Eq.~\eqref{eq:lr}. \algline{st:reszo:reg}
    \STATE Update $\bx_{t+1} = \bx_t - \eta\,\bg_t$. \algline{st:reszo:update}
\ENDFOR  \algline{st:reszo:end}
\end{algorithmic}
\end{algorithm}

\begin{remark}\label{remark:pert}
 (\textit{Adaptive Smoothing Radius}). It is noted that function evaluations are performed on the perturbed values ${\hat{\bm{x}}_{k}} = \bm{x}_{k} \!+\!\delta_k\bu_{k}$. Intuitively, adding random perturbations $\delta_k\bu_{k}$ promotes exploration and is critical for ensuring the algorithmic stability of the ReSZO method. The magnitude of perturbations is controlled by the smoothing radius parameter $\delta_k$. In general, a smaller $\delta$ reduces the bias in gradient estimation and improves convergence precision, whereas a larger $\delta$ promotes greater exploration during gradient estimation. To this end, in Step \ref{st:reszo:delta}, we propose an \textbf{adaptive perturbation scheme} \eqref{eq:adaptive} to effectively balance these two aspects based on the current optimization state. 
The key idea is to adjust the smoothing radius $\delta_t$ in proportion to the magnitude of $\|\bg_t\|$, given that the update in Step~\ref{st:reszo:update} follows $\| \bm{x}_{t}\!-\!\bm{x}_{t-1}\| = \eta \|\bg_t\|$. 
 Using a fixed $\delta$ would cause convergence to stall at an error level on the order of $\sO(\delta)$. In contrast, the adaptive scheme gradually reduces the perturbation size as the gradient norm decreases, enabling high convergence precision. See the proof of Lemma \ref{lemma:xi_smaller_than_grad_per_step} in the Appendix for more details, where this adaptive scheme is used to control the error magnitude of the gradient estimator. The necessity of introducing perturbations and the effectiveness of this adaptive scheme are further demonstrated through numerical experiments in Section~\ref{sec:pert}. 
\end{remark}

\subsection{Quadratic Regression-Based Single-Point Zeroth-Order Optimization (Q-ReSZO)}

In addition to fitting a linear model, other more sophisticated surrogate models, such as a \emph{quadratic} function, can be constructed based on available historical function evaluations to further enhance gradient estimation. At each iteration $t$, based on the second-order Taylor expansion at the current perturbed point $\hat{\bx}_t$, the available function evaluations $\{f(\hat{\bx}_k)\}_{k=0}^{t}$ can be expressed as: for all $ k \!=\! 0,1,\cdots,t$,
\begin{align} \label{eq:2taylor}
     f(\hat{\bx}_k) =  f(\hat{\bx}_t)\! +\! \nabla f(\hat{\bx}_t)^\top \big(\hat{\bx}_k \!-\! \hat{\bx}_t\big)  +  \frac{1}{2} \big(\hat{\bx}_k \!-\! \hat{\bx}_t\big)^\top H_f(\hat{\bx}_t) \big(\hat{\bx}_k \!-\! \hat{\bx}_t\big)
  \!+\! \epsilon_{t,k}(\hat{\bx}_k \!-\! \hat{\bx}_t), 
\end{align}
where $H_f(\cdot)$ is the Hessian matrix of function $f$, and
$\epsilon_{t,k}(\cdot)$ denotes the corresponding higher-order residual term. Accordingly, we can fit a quadratic surrogate function $f^\mathrm{s2}_t(\bx)$ 
\eqref{eq:sur2} to locally approximate the original function $f$ in \eqref{eq:problem}: 
\begin{align}\label{eq:sur2}
     \ f^\mathrm{s2}_t(\bx) \coloneq  f(\hat{\bx}_t) +  \bg_t^\top \big({\bx} - \hat{\bx}_t \big)  +  \frac{1}{2} \big({\bx} - \hat{\bx}_t \big)^\top \text{diag}(\bh_t) \big({\bx} - \hat{\bx}_t \big).  
\end{align}
Here, we use
the diagonal matrix $\text{diag}(\bm{h}_t)$ with vector $\bh_t\in\R^d$ to approximate the Hessian matrix, which reduces the number of variables for fitting. Similarly, 
a constant intercept term is omitted from the quadratic surrogate function \eqref{eq:sur2} to enforce $f^\mathrm{s2}_t(\hat{\bx}_t) = f(\hat{\bx}_t)$ and to mitigate model bias.
In addition,
 the recursive linear regression in a sliding window is employed, which uses only the latest $m$ function evaluations $\{\hat{\bm{x}}_{t-i+1}, f(\hat{\bm{x}}_{t-i+1})\}_{i=1}^{m}$  to fit the quadratic surrogate model \eqref{eq:sur2} at each iteration $t$. 
 
Then, we solve the following least squares problem \eqref{eq:ls:qr} to obtain the coefficients $\{\bg_t, \bh_t\}$:
\begin{align} \label{eq:ls:qr}
      \min_{\bm{g}_t,\bh_t}\, \sum_{i=1}^{m}  \Big[ f^\mathrm{s2}_t(\hat{\bx}_{t-i+1}) - f(\hat{\bx}_{t-i+1})\Big]^2  \    \Longrightarrow   \ 
    \begin{bmatrix}
      \bm{g}_t \\ \bh_t
    \end{bmatrix} = (\tilde{X}_t^\top \tilde{X}_t)^{\dagger}\tilde{X}_t^\top  {\bm{y}}_t,
\end{align}
where $\odot$ denotes the element-wise multiplication of two column vectors, $\tilde{X}_t \in \R^{(m-1) \times 2d}$ and ${\bm{y}}_t \in \R^{m-1}$ with
\begin{align}
\tilde{X}_t \!\coloneq\!\! \begin{bmatrix}
        {\bm{\Delta}^x_{t,1}}^\top &  \frac{1}{2}\big({\bm{\Delta}^x_{t,1}}\!\odot\!\bm{\Delta}^x_{t,1}\big)^\top\\  {\bm{\Delta}^x_{t,2}}^\top &  \frac{1}{2}\big({\bm{\Delta}^x_{t,2}}\!\odot\!\bm{\Delta}^x_{t,2}\big)^\top\\ \vdots &\vdots \\  {\bm{\Delta}^x_{t,m\!-\!1}}^\top&  \frac{1}{2}\big({\bm{\Delta}^x_{t,m\!-\!1}}\!\odot\!\bm{\Delta}^x_{t,m\!-\!1}\big)^\top 
    \end{bmatrix}\!,  {\bm{y}}_t \! \coloneq \!\! \begin{bmatrix}
       \Delta^f_{t,1}\\\Delta^f_{t,2}\\ \vdots\\ \Delta^f_{t,m\!-\!1}
    \end{bmatrix}\!.
\end{align}
Accordingly, our proposed \textbf{Q}uadratic \textbf{Re}gression-based \textbf{S}ingle-point \textbf{Z}eroth-order \textbf{O}ptimization (\textbf{Q-ReSZO}) method for solving problem \eqref{eq:problem} follows the same framework as the L-ReSZO method (i.e., Algorithm \ref{alg:lReSZO}), except that it solves for $\{\bg_t, \bh_t \}$ via \eqref{eq:ls:qr} in Step \ref{st:reszo:reg} and performs the iterative update by \eqref{eq:2update} with the gradient estimator $\nabla f^\mathrm{s2}_t(\bx)$ in Step \ref{st:reszo:update}. 
\begin{align} \label{eq:2update}
    \bx_{t+1} \!=\! \bx_t\! -\!\eta \nabla f^\mathrm{s2}_t(\bx_t) \!=\! \bx_t \!-\!\eta
   \big(\bg_t -  \delta_t \bh_t\odot \bu_t\big).
\end{align}





The ReSZO algorithms presented above aim for solving unconstrained optimization problems as formulated in \eqref{eq:problem};
nevertheless, they can be readily extended to constrained optimization problems using projection techniques \cite{chen2025continuous} and primal-dual algorithms \cite{chen2022model}.

\begin{remark}[\emph{Computational Cost versus Function-Evaluation Complexity}]
ReSZO requires computing a matrix inverse or pseudoinverse in the regression step at each iteration, introducing additional numerical computation compared with classical SZO and TZO methods. However, the central objective of ReSZO is to improve \emph{function-evaluation complexity}, which is often the dominant cost in black-box optimization. 
In many practical applications, such as real-world experiments, clinical studies, and large-scale or long-horizon simulations, obtaining a function evaluation can be substantially more expensive and time-consuming than performing an optimization update. For example, in power-system planning, evaluating a candidate solution may require extensive simulations across many operating conditions and contingencies, whereas the associated regression and update computations are comparatively inexpensive. 
ReSZO is particularly motivated by such settings: it reuses historical function evaluations to improve gradient estimation while requiring only one new function evaluation per iteration. As demonstrated in Section~\ref{sec:numerical}, ReSZO converges approximately two to three times faster than two-point ZO in terms of function queries, making the additional regression computation a favorable tradeoff when function-evaluation cost dominates.
\end{remark}

\section{Theoretical Convergence Analysis}\label{sec:analysis}

This section presents the non-asymptotic convergence analysis for the L-ReSZO method. Throughout the analysis, we use the notation $[N] \!\coloneqq\!\{1,2,\dots,N\}$, and $\xi_t \!\coloneqq\! \bm{g}_t \!-\! \nabla f(\bx_t)$  denotes the gradient estimation error at iteration $t$.

Throughout our analysis, we make the following standard assumption on the smoothness of the function $f$.

\begin{assumption}
\label{assumption:f_regularity}
The function $f$ is $L$-smooth, i.e., 
$
\|\nabla f(\bx)\!-\!\nabla f(\by)\|\leq L\|\bx-\by\|, \, \forall \bx,\by\in\mathbb{R}^d.
$
\end{assumption}

We recall the approximate gradient $\bm{g}_t$ can be written as equation \eqref{eqn:gt_in_analysis_defn}:
\begin{align}
\label{eqn:gt_in_analysis_defn}
    & \qquad\qquad\qquad\qquad \bm{g}_t = (S_t S_t^\top)^\dagger S_t \Delta_t^f, \\
    & S_t^\top \!:=\!\! \begin{bmatrix}
        (\hat{\bx}_{t-1} - \hat{\bx}_t)^\top \\
        (\hat{\bx}_{t-2} - \hat{\bx}_t)^\top \\
        \vdots \\
        (\hat{\bx}_{t-m+1}\!-\! \hat{\bx}_t)^\top \\
    \end{bmatrix}\!, \Delta_t^f\! := \!\!\begin{bmatrix}
        f(\hat{\bx}_{t-1}) \!-\! f(\hat{\bx}_t) \\
        f(\hat{\bx}_{t-2}) \!-\! f(\hat{\bx}_t) \\
        \vdots \\
        f(\hat{\bx}_{t-m+1}) \!-\! f(\hat{\bx}_t) \\
    \end{bmatrix}, \nonumber
\end{align}
where $S_t^\top =X_t \in \mathbb{R}^{(m-1) \times d }, \Delta_t^f=\bm{y}_t \in \mathbb{R}^{m-1}$. 


%
We also make the following assumption on the regression procedure. 
\begin{assumption}    \label{assumption:key_analysis_condition}
Suppose that $\bm{g}_t$ is as defined in (\ref{eqn:gt_in_analysis_defn}), and that $m > d$. We assume that
for any $t \geq m$,
    \begin{align}
    \label{eqn:analysis_key_assumption}
        \|\bm{g}_t - \nabla f(\hat{\bx}_t) \| \leq C_{d,m} \frac{L}{2} \max_{i=1:m-1}\|\hat{\bx}_{t-i} - \hat{\bx}_t\|,
    \end{align}
    for some constant $C_{d,m} > 0$ depending on the dimension $d$ and $m$. In particular, when the regression sample size is picked as $m = \sO(d)$, we denote $C_d := C_{d,m}$.
\end{assumption}

Assumption \ref{assumption:key_analysis_condition} bounds the error of the approximate gradient by the distance between points used in the regression $\hat{\bx}_{t-i}$ and the current point $\hat{\bx}_t$. For clarity, we note that the error $\bm{g}_t - \nabla f(\hat{\bx}_t)$ that appears in \eqref{eqn:analysis_key_assumption} is distinct from $\xi_t := \bm{g}_t - \nabla f(\bx_t)$, since the former is the bias compared to the gradient at the perturbed point while the latter is the bias compared to the gradient at the actual iterate of the algorithm.
We note that this assumption is difficult to rigorously prove, since it relies on the intricate behavior of the pseudoinverse of the random matrix $S_tS_t^\top$. We provide heuristic reasoning for its validity in Appendix \ref{app:assumption_2}, showing that $C_{d,m}$ scales at worst as $\sO(\sqrt{d}{\sqrt{m}}) =\sO(d)$ (when we set $m =\sO(d)$). In practice, we set $m = \sO(d)$ in our algorithm and present empirical evidence in Appendix C suggesting that in this setting, the factor $C_d := C_{d,m}$ in Assumption \ref{assumption:key_analysis_condition} actually scales on the order of $\sO(\sqrt{d})$. 


Before deriving our theoretical convergence results, we first specify the following condition on the initial $\tilde{\eta}$, namely the stepsize during the initial phase (before regression), indicating that it should be sufficiently small in a certain sense. Our subsequent results  require that $\tilde{\eta}$ satisfies this condition.

\begin{condition}
\label{condition:eta_tilde}
Recall that $\tilde{\eta}$ is the step-size of ReSZO during the initial update phase (before regression). We pick $\tilde{\eta}$ such that for any $t \in \{0,1,\dots,m-1\}$, it satisfies\footnote{We note that this is a relatively mild condition and expect it to hold for any sufficiently small $\tilde{\eta}$.} 
\begin{align}
\label{eq:eta_tilde_assumption}
    \|\bx_{t+1} - \bx_t\| = \|\tilde{\eta} \bm{g}_t\| \leq 2\eta \| \nabla f(\bx_t)\|,
\end{align}
where $\bm{g}_t$ denotes the estimator of $\nabla f(\bx_t)$ computed using the RSZO method in \eqref{eq:rfszo}.
\end{condition}

Our analysis is underpinned by the following key technical result, which bounds the gradient approximation error in terms of the norm of the true gradient. To be precise, we show that for any $0 < c < 1$, we can inductively guarantee that the gradient estimation bias $\|\xi_t \|$ always satisfies $\|\xi_t\| \leq c \|\nabla f(\bx_t)\|$ for all $t \geq m$. 

\begin{restatable}{lemma}{lemmaXiSmallerThanGrad}
\label{lemma:xi_smaller_than_grad_per_step}
Suppose that Assumptions \ref{assumption:f_regularity} and \ref{assumption:key_analysis_condition} hold, and that $m > d$. Suppose that
$\tilde{\eta}$ is picked such that Condition \ref{condition:eta_tilde} is satisfied. Then, for any $0 < c < 1$, by choosing $\eta = C\frac{1}{m C_d L}$ for some sufficiently small absolute constant $C > 0$, we have that 
\begin{align}
    \|\xi_t\| = \|\bm{g}_t - \nabla f(\bx_t)\| \leq c \|\nabla f(\bx_t) \|, 
\end{align}
for any $t \geq m$.
\end{restatable}
We defer the proof to Appendix \ref{app:Proofs}.

With Lemma \ref{lemma:xi_smaller_than_grad_per_step} in hand, we can guarantee the following function value improvement per step; the proof of this result is provided in Appendix \ref{app:Proofs}. 

\begin{restatable}{lemma}{fnImprovePerStep}
    \label{lemma:fn_improvement_per_step}
    Suppose that Assumptions \ref{assumption:f_regularity} and \ref{assumption:key_analysis_condition} hold, and that $m > d$.
Suppose that $\tilde{\eta}$ is picked such that Condition \ref{condition:eta_tilde} is satisfied. Then, by choosing $\eta = C\frac{1}{m C_d L}$ for some sufficiently small absolute constant $C > 0$, we have that 
\begin{align}
            f(\bx_{t+1}) - f(\bx_t) \leq - \frac{\eta}{8} \|\nabla f(\bx_t) \|^2, 
\end{align}
for any $t \geq m$.
\end{restatable}

With Lemma \ref{lemma:fn_improvement_per_step} in place, we can show by standard arguments that the L-ReSZO algorithm converges at a sublinear rate for smooth nonconvex functions and at a linear rate for smooth strongly convex functions, as stated in Theorems \ref{theorem:LReSZO-nonconvex-convergence} and \ref{theorem:LReSZO-strong-convex-convergence}. We defer the detailed proof to Appendix \ref{app:Proofs}.

\begin{restatable}{theorem}{thmNonconvex}
\label{theorem:LReSZO-nonconvex-convergence}
    Suppose that Assumptions \ref{assumption:f_regularity} and \ref{assumption:key_analysis_condition} hold, and that $d < m = \sO(d)$.
Suppose that $\tilde{\eta}$ is picked such that  Condition \ref{condition:eta_tilde} is satisfied. Then, by choosing $\eta = C\frac{1}{m C_d L}$ for some sufficiently small absolute constant $C > 0$, for any $T > m$,  the L-ReSZO algorithm achieves
\begin{align}
        \frac{1}{T}\sum_{t=m}^{T-1} \|\nabla f(\bx_{t})||^2 \leq \sO\left(\frac{d C_d L (f(\bx_{m}) - f^*)}{T} \right).
    \end{align}
\end{restatable}

This convergence rate scales as $\sO(\frac{d C_dL}{T})$ when we pick $m = \Theta(d)$, decaying at the rate $1/T$. We briefly comment on picking $m$ here. Our theory suggests that $m$ should be larger than $d$ to ensure that the matrix $S_t$ is full-rank, a requirement for our analysis of Assumption \ref{assumption:key_analysis_condition}); while it should  not exceed $\Theta(d)$, since $\eta$ scales as $\sO(1/m)$ and increasing $m$ slows down the convergence rate.

Next, we present the linear convergence result for strongly convex functions.

\begin{restatable}{theorem}{thmStrongConvex}
\label{theorem:LReSZO-strong-convex-convergence}
    Suppose that Assumptions \ref{assumption:f_regularity} and \ref{assumption:key_analysis_condition} hold, and that $d < m = \sO(d)$.
Suppose that $\tilde{\eta}$ is picked such that Condition \ref{condition:eta_tilde} is satisfied. Then, by choosing $\eta = C\frac{1}{m C_d L}$ for some sufficiently small absolute constant $C > 0$, for any $T > m$, the L-ReSZO algorithm achieves
        \begin{align}
        f(\bx_T)  - f^* \leq \left(1  - \sO\big(\frac{\mu}{C_dL d} \big) \right)^{T-m} (f(\bx_m) - f^*).
    \end{align}
\end{restatable}

Theorem \ref{theorem:LReSZO-strong-convex-convergence} shows that the L-ReSZO algorithm achieves a linear convergence rate for strongly convex functions, which is slower than that of two-point ZO methods only by a factor of $C_d$. 


\section{Numerical Experiments} \label{sec:numerical}

This section demonstrates the performance of the proposed ReSZO methods in comparison with classic ZO methods through numerical experiments across  various optimization problems, 
and studies the impact of adding perturbations.

\subsection{Comparison with Classic ZO Methods} \label{sec:testcom}

We compare the performance of the proposed L-ReSZO (i.e., Algorithm \ref{alg:lReSZO}), Q-ReSZO \eqref{eq:2update} methods with the RSZO \eqref{eq:rfszo} and TZO \eqref{eq:tzo} methods 
across four optimization tasks, including 
ridge regression, logistic regression, Rosenbrock function minimization, and neural networks training. 

\vspace{3pt}
\noindent\textbf{Case (a) Ridge Regression}.
Consider solving the ridge regression problem \eqref{eq:ridge} \cite{uribe2020dual}:
\begin{align}\label{eq:ridge}
    \min_{\bx\in\R^d}\, f(\bx)= \frac{1}{2}||\bm{y} - H\bx ||_2^2 + \frac{\lambda}{2}||\bx||_2^2.
\end{align}
In our experiments,
each entry of the matrix $H\in \R^{N\times d}$ is independently drawn from the standard Gaussian distribution $\mathcal{N}(0,1)$. The column vector $\bm{y}\in \R^N$ is constructed by letting
$\bm{y}= \frac{1}{2}H\bm{1}_d +\bm{\epsilon}$, where $\bm{1}_d$ denotes the $d$-dimensional all-ones vector and $\bm{\epsilon}\sim \mathcal{N}(\bm{0},0.1\bm{I})$ represents the additive Gaussian noise.  The regularization parameter is fixed as $\lambda=0.1$, and the initial point is set to $\bx_0 = \bm{0}_d$. We set 
 $N = 1500$ and the problem dimension $d=500$.

\vspace{3pt}
\noindent
 \textbf{Case (b) Logistic Regression}. 
 Consider solving the logistic regression  problem \eqref{eq:logis} \cite{uribe2020dual}:
\begin{align} \label{eq:logis}
    \min_{\bx\in\R^d} f(\bx) = \frac{1}{2} \sum_{i=1}^N \log\big(  1+\exp(-y_i\cdot \bs_i^\top \bx) \big) + \frac{\lambda}{2}||\bx||_2^2,
\end{align}
where  $\bs_i\in \R^d$ is one of the data samples, and $y_i\in\{-1,1\}$ is the corresponding class label. Each element of a data sample $\bs_i$ is independently drawn from the uniform distribution $\mathrm{Unif}([-1,1])$, and the corresponding label is computed as $y_i=\mathrm{sign}(\frac{1}{2}\bs_i^\top \mathbf{1}_d)$. The regularization parameter is fixed as $\lambda=0.1$. We set    $N=1000$ and $d =100$.

\vspace{3pt}
\noindent
\textbf{Case (c) Rosenbrock Function}. Consider minimizing the nonconvex Rosenbrock function \eqref{eq:rosen}: 
\begin{align}\label{eq:rosen}
    \min_{\bx\in\R^d} f(\bx) \!=\! \sum_{i=1}^d\! \big[ 100\big((x_i\!+\!1)^2\!-\!x_{i+1}-1\big)^2 \!+\! x_i^2 \big],
\end{align}
which is a well-known benchmark problem for numerical optimization \cite{shang2006note}. The global optimal solution and objective value are $\bx^*\!=\!\bm{0}$ and $f^* \!=\! 0$. We set $d\!=\!200$ and the initial point $\bx_0 \!=\! 0.5\bm{1}_d$. 

\vspace{3pt}
\noindent
\textbf{Case (d) Neural Network Training}. 
Consider training a three-layer fully connected neural network by minimizing the loss function \eqref{eq:nn}:
\begin{align}\label{eq:nn}
  \min_{\bx\in\R^d}  f(\bx) = \sum_{i=1}^N \Big[ \bm{w}_o^\top \sigma(W_3\sigma (W_2\sigma(W_1\bs_i + \bm{b}_1)  +\bm{b}_2) +\bm{b}_3)   - y_i  \Big]^2, 
\end{align}
where $\bx\coloneq (W_1, W_2, W_3, \bm{b}_1, \bm{b}_2, \bm{b}_3,\bm{w}_o)$, and $\sigma(\cdot)$ denotes the Sigmoid activation function. Let $\bs_i\in\R^n$ and $y_i\in\R$ be the input and output of each sample $i$ with $n = 6$. Thus, the problem dimension is $d = 3n^2+4n = 132$. We set the sample size to $N=500$. Each sample input $\bs_i$ is randomly drawn from the standard Gaussian distribution $\mathcal{N}(0,1)$, and the corresponding output is computed by the neural network with ground truth $\bx^*$, where each entry of $\bx^*$ is randomly generated from $\mathcal{N}(0,1)$. Let the initial point be $\bx_0 \!=\! \bx^* + \bm{\varepsilon}$, and each entry of $\bm{\varepsilon}$ is randomly drawn from the uniform distribution $\mathrm{Unif}([-1,1])$.

\begin{figure*}[t]
    \centering
    \includegraphics[width=\textwidth]{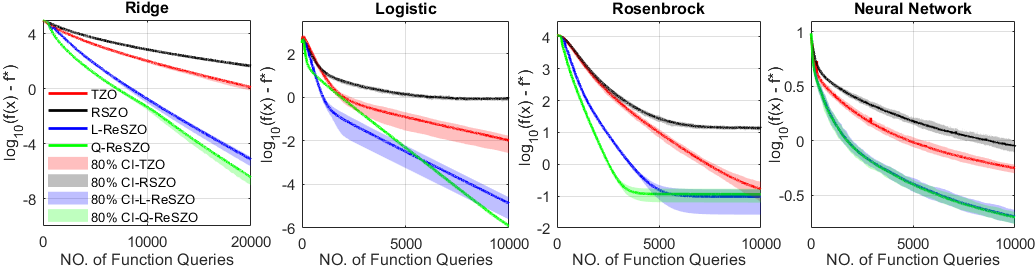}
    \caption{The convergence and 80\% confidence intervals (CI) of TZO \eqref{eq:tzo}, RSZO \eqref{eq:rfszo}, L-ReSZO (Algorithm \ref{alg:lReSZO}) and Q-ReSZO \eqref{eq:2update} for solving ridge regression \eqref{eq:ridge}, logistic regression \eqref{eq:logis}, Rosenbrock function minimization \eqref{eq:rosen}, and neural network training \eqref{eq:nn}. }
    \label{fig:combine}
\end{figure*}

We optimize the stepsize $\eta$ for each ZO method and the smoothing radius $\delta$ for RSZO and TZO via grid search to achieve the fastest convergence. The adaptive smoothing radius scheme \eqref{eq:adaptive} is employed in L-ReSZO and Q-ReSZO. The selected parameters of $\eta$ and  $\delta$, as well as the sample size $m$ for regression, are presented in Table~\ref{tab:para} in Appendix \ref{app:setting}. 
For each of the four optimization tasks described above, we run each ZO method 50 times and calculate the mean and $80\%$-confidence interval (CI) of the optimality gap $f(\bx) \!-\! f^*$. Figure \ref{fig:combine} illustrates the convergence results using the TZO, RSZO, L-ReSZO, and Q-ReSZO methods, evaluated based on the number of function queries.

As shown in Figure~\ref{fig:combine}, the proposed L-ReSZO and Q-ReSZO methods converge two to three times faster than TZO in terms of function query complexity and significantly outperform RSZO across all test cases. With respect to iteration complexity, L-ReSZO and Q-ReSZO achieve convergence rates comparable to, or slightly faster than, TZO. Since L-ReSZO and Q-ReSZO require only one function query per iteration, they are approximately two to three times more efficient than TZO in terms of function query complexity. In addition, Q-ReSZO exhibits slightly faster convergence than L-ReSZO, which can be attributed to its more accurate gradient estimation enabled by quadratic surrogate functions.


\subsection{Effect of Perturbations on ReSZO}\label{sec:pert}

To analyze the effect of perturbations used in function evaluations, we evaluate the proposed L-ReSZO and Q-ReSZO methods under varying levels of the smoothing radius 
$\delta$ when solving the ridge regression problem~\eqref{eq:ridge}. The experimental setup remains the same as that described in the previous subsection. Specifically, $\delta$ is progressively reduced from 1, 0.1, 0.01, 0.001, to zero (no perturbation), and
we also perform the adaptive smoothing radius scheme \eqref{eq:adaptive}.

\begin{figure}
    \centering
        \includegraphics[width=0.6\textwidth]{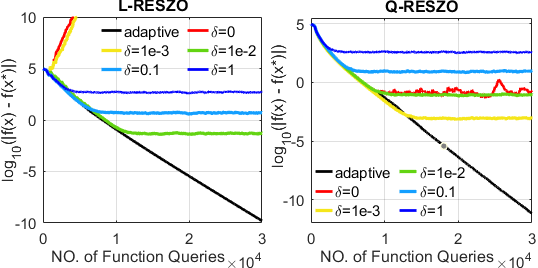}
    \caption{Convergence of L-ReSZO and Q-ReSZO for solving the ridge regression problem \eqref{eq:ridge}
    under different levels of perturbations and using the adaptive smoothing radius scheme \eqref{eq:adaptive}. }
    \label{fig:per}
    \vspace{-10pt}
\end{figure}

The test results are illustrated in Figure~\ref{fig:per}. It is observed that a smaller smoothing radius 
$\delta$ generally leads to improved convergence precision. This is because ZO methods estimate the gradient at a perturbed point, i.e.,  $\nabla f(\hat{\bx}) = \nabla f(\bx + \delta \bu_t)$, and thus a smaller $\delta$ results in lower gradient estimation errors. However, the introduction of random perturbations is necessary for the stability of ReSZO. As shown in Figure~\ref{fig:per}, when no perturbation is applied ($\delta=0$) or when the perturbation is excessively small ($\delta=0.001$), the L-ReSZO method diverges, and the Q-ReSZO method exhibits large oscillations in the former case. The proposed adaptive smoothing radius scheme effectively addresses these issues by automatically adjusting $\delta$, achieving both algorithmic stability and high convergence precision.

\section{Conclusion}

This paper proposes a novel yet simple single-point zeroth-order optimization (SZO) framework termed Regression-based SZO (ReSZO), which significantly improves gradient estimation by leveraging historical function evaluations from previous iterations.
We theoretically establish that the ReSZO method achieves convergence rates comparable to those of the two-point ZO method for both smooth nonconvex and strongly convex functions. Furthermore, extensive numerical experiments demonstrate that ReSZO not only matches but empirically converges approximately two to three times faster than two-point ZO methods in terms of function query complexity across all test cases. To  our knowledge, this is the first single-point zeroth-order method that attains convergence rates comparable to two-point methods while achieving superior empirical performance. Future work will explore strategies to further enhance the ReSZO framework, including more expressive surrogate models, weighted regression schemes, and extensions to noisy function evaluations that better reflect practical settings.

\section*{Acknowledgments}

This work was supported in part by NSF AMPS 2523934, and in part by NSF CAREER 2541998.

\bibliographystyle{unsrt}

\bibliography{reference}

\newpage

\appendix

\section{Numerical Experiment Settings}\label{app:setting}

For the numerical experiments in Section \ref{sec:testcom}, 
we optimize the stepsize $\eta$ for each ZO method and the smoothing radius $\delta$ for RSZO and TZO via grid search to achieve the fastest convergence. The adaptive smoothing radius scheme \eqref{eq:adaptive} is employed in L-ReSZO and Q-ReSZO. The selected parameters of $\eta$ and  $\delta$, as well as the sample size $m$ for regression, are presented in Table~\ref{tab:para}. 

The same sample size $m$ is used for both L-ReSZO and Q-ReSZO, with $m$ set to 510, 110, 210, and 6 for Cases (a)–(d), respectively. In Cases (a)–(c), the sample size is chosen as $m=d+10$, slightly exceeding the problem dimension $d$. In Case (d), a smaller sample size yields better convergence, and thus 
$m=6$ is adopted. This observation suggests that, for certain problems, only a small number of recent function evaluations is sufficient to construct an accurate gradient estimate using the proposed ReSZO method.

\begin{table}[H]
\centering
  \caption{The selected smoothing radius $\delta$, stepsize $\eta$, and sample size $m$ for the ZO methods in Cases (a)-(d) in Section \ref{sec:testcom}.} \label{tab:para}
\resizebox{\textwidth}{!}{%
\begin{tabular}{l ccc ccc ccc ccc}
  \toprule[1.2pt]
  \multirow{2}{*}{Methods} &
  \multicolumn{3}{c}{Ridge Reg.} &
  \multicolumn{3}{c}{Logistic Reg.} &
  \multicolumn{3}{c}{Rosenbrock Func.} &
  \multicolumn{3}{c}{Neural Network} \\
  \cmidrule(r){2-4}\cmidrule(r){5-7}\cmidrule(r){8-10}\cmidrule(r){11-13}
  & $\boldsymbol{\delta}$ & $\boldsymbol{\eta}$ & $\boldsymbol{m}$
  & $\boldsymbol{\delta}$ & $\boldsymbol{\eta}$ & $\boldsymbol{m}$
  & $\boldsymbol{\delta}$ & $\boldsymbol{\eta}$ & $\boldsymbol{m}$
  & $\boldsymbol{\delta}$ & $\boldsymbol{\eta}$ & $\boldsymbol{m}$ \\
  \midrule
  TZO  & 0.01 & $1.2\times10^{-6}$ & -- & 0.01 & $1.6\times10^{-3}$ & -- & 0.01 & $4.5\times10^{-6}$ & -- & 0.01 & $3.8\times10^{-4}$ & -- \\
  RSZO & 0.2   & $3\times10^{-7}$ & -- & 2    & $5\times10^{-4}$   & -- & 0.5  & $2\times10^{-6}$   & -- & 0.05 & $1.1\times10^{-4}$ & -- \\
  L-ReSZO & adaptive & $1.5\times10^{-6}$ & 510 & adaptive  & $1.8\times10^{-3}$ & 110 & adaptive  & $4.2\times10^{-6}$ & 210 & adaptive  & $1.7\times10^{-3}$ & 6 \\
  Q-ReSZO & adaptive  & $2\times10^{-6}$ & 510 & adaptive  & $5\times10^{-3}$ & 110 & adaptive  & $1\times10^{-5}$ & 210 & adaptive  & $1.8\times10^{-3}$ & 6 \\
  \bottomrule[1.2pt]
\end{tabular}%
}
\end{table}

 \section{Proofs for Theoretical Convergence Analysis} \label{app:Proofs}

 \subsection{Theoretical justification for Assumption \ref{assumption:key_analysis_condition}}
 \label{app:assumption_2}

To facilitate our analysis, we first show the following result which provides an alternative characterization of the approximation error term $\|\bm{g}_t - \nabla f(\hat{\bx}_t)\|$ that appears in Assumption \ref{assumption:key_analysis_condition}. 
\begin{lemma}
    \label{lemma:alt_form_gt_grad_diff}
    Suppose $\bm{g}_t$ is defined as in \eqref{eqn:gt_in_analysis_defn}, i.e. 
\begin{align*}
    & \ \bm{g}_t = (S_t S_t^\top)^\dagger S_t \Delta_t^f, \quad S_t^\top := \begin{bmatrix}
        (\hat{\bx}_{t-1} - \hat{\bx}_t)^\top \\
        (\hat{\bx}_{t-2} - \hat{\bx}_t)^\top \\
        \vdots \\
        (\hat{\bx}_{t-m+1} - \hat{\bx}_t)^\top \\
    \end{bmatrix},
     \quad \Delta_t^f := \begin{bmatrix}
        f(\hat{\bx}_{t-1}) - f(\hat{\bx}_t) \\
        f(\hat{\bx}_{t-2}) - f(\hat{\bx}_t) \\
        \vdots \\
        f(\hat{\bx}_{t-m+1}) - f(\hat{\bx}_t) \\
    \end{bmatrix}.
\end{align*}
Define the errors of the first-order Taylor expansion as:
\begin{align}
    \label{eq:epsilon_t_defn}
\epsilon_t := \begin{bmatrix}
        f(\hat{\bx}_{t-1}) - (f(\hat{\bx}_t) + \langle \nabla f(\hat{\bx}_t), \hat{\bx}_{t-1} - \hat{\bx}_t\rangle) \\
        \vdots \\
        f(\hat{\bx}_{t-m+1}) - (f(\hat{\bx}_t) + \langle \nabla f(\hat{\bx}_t), \hat{\bx}_{t-m+1} - \hat{\bx}_t\rangle)
    \end{bmatrix}.
\end{align}
Suppose $S_tS_t^\top \in \mathbb{R}^{d \times d}$ is invertible.  
Then,
\begin{align*}
    \bm{g}_t - \nabla f(\hat{\bx}_t) = (S_t S_t^\top)^{-1} S_t \epsilon_t.
\end{align*}
\end{lemma}
\begin{proof}
    Using the assumption that $S_t S_t^\top \in \R^{d \times d}$ is invertible, we have
    \begin{align*}
        \bm{g}_t = & \ (S_t S_t^\top)^{-1}S_t \Delta_t^f = \ (S_t S_t^\top)^{-1}S_t \left(S_t^\top \nabla f(\hat{\bx}_t) + \epsilon_t \right) \\
        = & \ \nabla f(\hat{\bx}_t) + (S_t S_t^\top)^{-1} S_t \epsilon_t,
    \end{align*}
    as desired.
\end{proof}

We can further simplify the above result to an even simpler characterization.

\begin{lemma}
    \label{lemma:alt_form_gt_grad_diff_v2}
    Consider the cumulative sum matrix $C \in \mathbb{R}^{m \times m}$, defined as
\begin{align*}
    C \;=\;
\begin{pmatrix}
1 & 0 & 0 & \cdots & 0 \\
1 & 1 & 0 & \cdots & 0 \\
1 & 1 & 1 & \cdots & 0 \\
\vdots & \vdots & \vdots & \ddots & \vdots \\
1 & 1 & 1 & \cdots & 1
\end{pmatrix}.
\end{align*}
Define also the matrix 
\begin{align*}
D_t^\top := \begin{bmatrix}
        (\hat{\bx}_{t-1} - \hat{\bx}_t)^\top \\
        (\hat{\bx}_{t-2} - \hat{\bx}_{t-1})^\top \\
        \vdots \\
        (\hat{\bx}_{t-m+1} - \hat{\bx}_{t-m+2})^\top \\
    \end{bmatrix}, \\
\end{align*}
as well as the following error vector
\begin{align}
    \label{eq:eps_hat_t_defn}
\hat{\epsilon}_t := \begin{bmatrix}
        f(\hat{\bx}_{t-1}) - (f(\hat{\bx}_t) + \langle \nabla f(\hat{\bx}_t), \hat{\bx}_{t-1} - \hat{\bx}_t\rangle)  \\
        f(\hat{\bx}_{t-2}) - \left(f(\hat{\bx}_{t-1}) + \langle \nabla f(\hat{\bx}_t), \hat{\bx}_{t-2} - \hat{\bx}_{t-1}\rangle\right)  \\
        \vdots \\
        f(\hat{\bx}_{t-m+1}) - \left(f(\hat{\bx}_{t-m+2}) + \langle \nabla f(\hat{\bx}_t), \hat{\bx}_{t-m+1} - \hat{\bx}_{t-m+2}\rangle\right)
    \end{bmatrix}.
\end{align}
Suppose $D$ has full rank. Then,
\begin{align*}
    \bm{g}_t - \nabla f(\hat{\bx}_t) = (D_tD_t^\top)^{-1} D_t\hat{\epsilon}_t.
\end{align*}
\end{lemma}
\begin{proof}
    From the earlier result, we saw that when $S$ is full rank (equivalent to $D$ being full rank), we have
        \begin{align*}
        \bm{g}_t = & \ (S_t S_t^\top)^{-1}S_t \Delta_t^f = \ (S_t S_t^\top)^{-1}S_t \left(S_t^\top \nabla f(\hat{\bx}_t) + \epsilon_t \right) \\
        = & \ \nabla f(\hat{\bx}_t) + (S_t S_t^\top)^{-1} S_t \epsilon_t.
    \end{align*}
    It can be easily verified that $S_t^\top = C D_t^\top$. Since $(S_t^\top)^\dagger = (CD_t^\top)^\dagger  = (D_t^\top)^\dagger C^{-1}$, we have that 
        \begin{align*}
        \bm{g}_t = & \ \nabla f(\hat{\bx}_t) + (S_t S_t^\top)^{-1} S_t \epsilon_t \\ 
        = & \ \nabla f(\hat{\bx}_t) + (D_t D_t^\top)^{-1} D_t C^{-1} \epsilon_t.
    \end{align*}
    Since $C^{-1}$ is the following matrix 
    \begin{align*}
        C^{-1} \;=\;
\begin{pmatrix}
1 & 0 & 0 & \cdots & 0 \\
-1 & 1 & 0 & \cdots & 0 \\
0 & -1 & 1 & \cdots & 0 \\
\vdots & \ddots & \ddots & \ddots & 0 \\
0 & \cdots & 0 & -1 & 1
\end{pmatrix},
    \end{align*}
    the desired result follows since $C^{-1} \epsilon_t = \hat{\epsilon}_t.$
\end{proof}

In order for Assumption \ref{assumption:key_analysis_condition} to hold, we need to rely on the following assumption on the smallest singular value of $D_t$. While difficult to prove since it involves computing the smallest singular value of random matrices with mild correlation across rows, we provide heuristic support for its validity.

\begin{assumption}
\label{assumption:Dt_smallest_singular_value}
    Let $D_t^\top \in \mathbb{R}^{(m-1) \times d}$ be the matrix 
\begin{align*}
D_t^\top := \begin{bmatrix}
        (\hat{\bx}_{t-1} - \hat{\bx}_t)^\top \\
        (\hat{\bx}_{t-2} - \hat{\bx}_{t-1})^\top \\
        \vdots \\
        (\hat{\bx}_{t-m+1} - \hat{\bx}_{t-m+2})^\top \\
    \end{bmatrix}. \\
\end{align*}
Then, the smallest singular value of $D_t^\top$ satisfies the following bound 
$$\sigma_{\min}(D_t^\top) \geq \Omega\left(\frac{1}{\Gamma_{d,m}} \sqrt{\sum_{i=1}^{m-1} \|\hat{\bx}_{t-i} - \hat{\bx}_{t-i+1} \|^2} \right).$$
Here $\Gamma_{d,m}$ is a dimension-dependent constant satisfying
\[
\Gamma_{d,m} = \Theta\!\left(\sqrt{d}\sqrt{m}\right).
\]
\end{assumption}

\begin{proof}[Theoretical justification for Assumption \ref{assumption:Dt_smallest_singular_value}]
We note that the matrix $D_t^\top$ can also be written as 
\begin{align*}
    D_t^\top = & \  \begin{bmatrix}
        (\hat{\bx}_{t-1} - \hat{\bx}_t)^\top \\
        (\hat{\bx}_{t-2} - \hat{\bx}_{t-1})^\top \\
        \vdots \\
        (\hat{\bx}_{t-m+1} - \hat{\bx}_{t-m+2})^\top \\
    \end{bmatrix} \\
    = & \ \begin{bmatrix}
        (\hat{\bx}_{t-1} - {\bx}_{t-1})^\top +  ({\bx}_{t-1} - {\bx}_{t})^\top + ({\bx}_{t} - \hat{\bx}_{t})^\top \\
        (\hat{\bx}_{t-2} - {\bx}_{t-2})^\top +  ({\bx}_{t-2} - {\bx}_{t-1})^\top + ({\bx}_{t-1} - \hat{\bx}_{t-1})^\top \\
        \vdots \\
        (\hat{\bx}_{t-m+1} - {\bx}_{t-m+1})^\top +  ({\bx}_{t-m+1} - {\bx}_{t-m+2})^\top + ({\bx}_{t-m+2} - \hat{\bx}_{t-m+2})^\top \\
    \end{bmatrix} \\
        = & \ \begin{bmatrix}
        (\bu_{t-1})^\top +  ({\bx}_{t-1} - {\bx}_{t})^\top + (\bu_{t})^\top \\
        (\bu_{t-2})^\top +  ({\bx}_{t-2} - {\bx}_{t-1})^\top + (\bu_{t-1})^\top \\
        \vdots \\
        (\bu_{t-m+1})^\top +  ({\bx}_{t-m+1} - {\bx}_{t-m+2})^\top + (\bu_{t-m+2})^\top \\
    \end{bmatrix} \\
    =  & \ G_t^\top + H_t^\top,
\end{align*}
where 
\begin{align*}
    G_t^\top = & \ \begin{bmatrix}
        ({\bx}_{t-1} - {\bx}_t)^\top \\
        ({\bx}_{t-2} - {\bx}_{t-1})^\top \\
        \vdots \\
        ({\bx}_{t-m+1} - {\bx}_{t-m+2})^\top
    \end{bmatrix}, \quad
    H_t^\top = \begin{bmatrix}
        (\bu_{t-1})^\top - (\bu_{t})^\top \\
        (\bu_{t-2})^\top -  (\bu_{t-1})^\top \\
        \vdots \\
        (\bu_{t-m+1})^\top - (\bu_{t-m+2})^\top \\
    \end{bmatrix},
\end{align*}
with $\bm{u_{j}}^\top$ is sampled independently and uniformly from the sphere with norm $\|\bm{x}_{j} - \bm{x}_{j-1} \|$. To find the smallest singular value of $D_t^\top$, we thus have to find the smallest eigenvalue of the matrix $(D_t D_t^\top) = (G_t G_t^\top + G_t H_t^\top + H_t G_t^\top + H_t H_t^\top)$. While $G_t$ and $H_t$ are slightly correlated, the cross-terms should display concentration effects and thus dominated by the effects of the terms $G_t G_t^\top + H_tH_t^\top $. To see this, taking $G_t H_t^\top$ for example, we observe that 
\begin{align}
    G_t H_t^\top = \sum_{j=1}^{m-1} (\bm{x}_{t-j} - \bm{x}_{t-j+1})(\bm{u}_{t-j} - \bm{u}_{t-j+1})^\top.
\end{align}
Since $\bu_{t-j+1}$ is independent of $\bm{x}_{t-j} - \bm{x}_{t-j+1}$, the term $\sum_{j=1}^{m-1} (\bm{x}_{t-j} - \bm{x}_{t-j+1})(\bm{u}_{t-j+1})^\top$ is mean zero and exhibits concentration of measure effects. Similarly, any correlation between $\bm{x}_{t-j} - \bm{x}_{t-j+1}$ and $\bm{u}_{t-j}$ is likely mild; when the iterates are generated before the regression starts, e.g., via one-point method, $\bm{x}_{t-j} - \bm{x}_{t-j+1}$ will have no relation to $\bu_{t-j}$, since the random direction at iterate $t-j$ is independent of the perturbation $\bu_{t-j}$ which is only to obtain the function value evaluation. When $\bm{x}_{t-j} - \bm{x}_{t-j+1} = \eta\bm{g}_{t-j}$ is the result of a regression-based gradient estimation update, $\bm{g}_{t-j}$ is affected by $\bm{u}_{t-j}$, but only in one row, and thus unlikely to have significant correlation. 

We may thus assume that 
$$ \sigma_{\min}(D_t D_t^\top) = \Theta( \sigma_{\min}(G_t G_t^\top + H_t H_t^\top)) \geq \Theta( \sigma_{\min}(H_t H_t^\top)).$$
Fully characterizing the distribution of $\sigma_{\min}(H_t H_t^\top)$ is highly non-trivial since the norm of each row of $H_t^\top$ has a different variance, namely that the $j$-th row of $H_t^\top$, which is equal to $(\bm{u}_{t-j} - \bm{u}_{t-j+1})^\top$, has variance $\|\bm{u}_{t-j} \|^2 + \|\bm{u}_{t-j+1}\|^2 = \|\bm{x}_{t-j} - \bm{x}_{t-j-1} \|^2  + \|\bm{x}_{t-j+1} - \bm{x}_{t-j} \|^2$. 

To proceed, note that $H_t^\top \in \mathbb{R}^{(m-1)\times d}$ has rows
$h_j^\top := (\bu_{t-j}-\bu_{t-j+1})^\top$ for $j=1,\dots,m-1$.
To build intuition for the scale of $\sigma_{\min}(H_t^\top)$, consider the Gram matrix
\[
H_t H_t^\top \;=\; \sum_{j=1}^{m-1} h_j h_j^\top
\;=\; \sum_{j=1}^{m-1} (\bu_{t-j}-\bu_{t-j+1})(\bu_{t-j}-\bu_{t-j+1})^\top .
\]
Expanding the summand gives
\[
(\bu_{t-j}-\bu_{t-j+1})(\bu_{t-j}-\bu_{t-j+1})^\top
=
\bu_{t-j}\bu_{t-j}^\top + \bu_{t-j+1}\bu_{t-j+1}^\top
-\bu_{t-j}\bu_{t-j+1}^\top - \bu_{t-j+1}\bu_{t-j}^\top .
\]
Summing over $j$ yields
\[
H_t H_t^\top
=
\underbrace{\sum_{j=1}^{m-1}\bu_{t-j}\bu_{t-j}^\top + \sum_{j=1}^{m-1}\bu_{t-j+1}\bu_{t-j+1}^\top}_{\text{``diagonal'' terms}}
\;-\;
\underbrace{\sum_{j=1}^{m-1}\big(\bu_{t-j}\bu_{t-j+1}^\top + \bu_{t-j+1}\bu_{t-j}^\top\big)}_{\text{cross terms}}.
\]
The diagonal part is essentially a weighted sum of rank-one matrices of the form $\bu_k\bu_k^\top$.
The cross terms are mean-zero (conditionally on $\{\bx_k\}$) because the directions $\bu_k$
are independent across $k$ and isotropic, so $\mathbb{E}[\bu_{t-j}\bu_{t-j+1}^\top]=0$.
Consequently, by concentration of measure, the cross terms exhibit cancellations and are typically of
smaller spectral norm compared to the diagonal terms. This motivates the approximation
\[
\sigma_{\min}(H_t H_t^\top)
\;\approx\;
\sigma_{\min}\!\left(\sum_{k=t-m+2}^{t-1} \bu_k \bu_k^\top\right).
\]

Next, assume that $\|\bu_k\|^2$ are comparable over the window, i.e.,
$\|\bu_k\|^2 \asymp r^2$ for $k=t-m+2,\dots,t-1$.
Under either the Gaussian model $\bu_k \sim \mathcal{N}(0,\tfrac{r^2}{d}I_d)$ or the uniform-sphere model
(with comparable sub-Gaussian tails), the matrix with rows proportional to $\bu_k^\top$ is a
(sub-)Gaussian design with unit-variance entries up to scaling.
Therefore, by classic results on the smallest singular value of sub-Gaussian random matrices
(e.g., Theorem~1.1 in \cite{rudelson2009smallest}), we have heuristically that
\[
\sigma_{\min}\!\left(\sum_{k=t-m+2}^{t-1} \bu_k \bu_k^\top\right)
\;=\;
\Theta\!\left(\frac{r^2}{d}\right),
\]
which implies
\[
\sigma_{\min}(H_t^\top)
\;=\;
\sqrt{\sigma_{\min}(H_t H_t^\top)}
\;\approx\;
\Theta\!\left(\frac{r}{\sqrt d}\right).
\]
Finally, using $\|\bu_k\|^2\asymp r^2$ and $m-1$ terms in the window,
\[
r \;\asymp\; \sqrt{\frac{1}{m-1}\sum_{k=t-m+2}^{t-1}\|\bu_k\|^2},
\]
so the above scaling can be rewritten as
\begin{align*}
\sigma_{\min}(H_t^\top)
\;\approx\; & \ 
\Theta\!\left(\frac{1}{\sqrt d}\sqrt{\frac{1}{m-1}\sum_{k=t-m+2}^{t-1}\|\bu_k\|^2}\right)
\;=\;
\Theta\!\left(\frac{1}{\sqrt{d(m-1)}}\sqrt{\sum_{k=t-m+2}^{t-1}\|\bu_k\|^2}\right) \\
= & \ \Omega\left(\frac{1}{\sqrt{d}\sqrt{m}} \sqrt{\sum_{i=1}^{m-1} \|\hat{\bx}_{t-i} - \hat{\bx}_{t-i+1} \|^2} \right).
\end{align*}


\end{proof}

We are now ready to prove Assumption \ref{assumption:key_analysis_condition}, assuming that Assumption \ref{assumption:Dt_smallest_singular_value} holds.

\begin{proof}[Proof of Assumption \ref{assumption:key_analysis_condition}]
    Suppose Assumption \ref{assumption:Dt_smallest_singular_value} holds. Recall from Lemma \ref{lemma:alt_form_gt_grad_diff_v2} that 
    \begin{align*}
    \bm{g}_t - \nabla f(\hat{\bx}_t) = (D_tD_t^\top)^{-1} D_t\hat{\epsilon}_t, 
\end{align*}
where 
\begin{align*}
D_t^\top := \begin{bmatrix}
        (\hat{\bx}_{t-1} - \hat{\bx}_t)^\top \\
        (\hat{\bx}_{t-2} - \hat{\bx}_{t-1})^\top \\
        \vdots \\
        (\hat{\bx}_{t-m+1} - \hat{\bx}_{t-m+2})^\top \\
    \end{bmatrix}, \\
\end{align*}
and
\begin{align*}
\hat{\epsilon}_t := \begin{bmatrix}
        f(\hat{\bx}_{t-1}) - (f(\hat{\bx}_t) + \langle \nabla f(\hat{\bx}_t), \hat{\bx}_{t-1} - \hat{\bx}_t\rangle)  \\
        f(\hat{\bx}_{t-2}) - \left(f(\hat{\bx}_{t-1}) + \langle \nabla f(\hat{\bx}_t), \hat{\bx}_{t-2} - \hat{\bx}_{t-1}\rangle\right)  \\
        \vdots \\
        f(\hat{\bx}_{t-m+1}) - \left(f(\hat{\bx}_{t-m+2}) + \langle \nabla f(\hat{\bx}_t), \hat{\bx}_{t-m+1} - \hat{\bx}_{t-m+2}\rangle\right)
    \end{bmatrix}.
\end{align*}
From Assumption \ref{assumption:Dt_smallest_singular_value}, we have that 
$$\sigma_{\min}(D_t^\top) \geq \Omega\left(\frac{1}{\Gamma_{d,m}} \sqrt{\sum_{i=1}^{m-1} \|\hat{\bx}_{t-i} - \hat{\bx}_{t-i+1} \|^2} \right).$$
It follows then that
\begin{align*}
    \|\bm{g}_t - \nabla f(\hat{x}_t)\| \leq \frac{1}{\sigma_{\min}(D_t)} \|\hat{\epsilon}_t \| \leq \frac{\sO(\sqrt{d}\sqrt{m})}{\sqrt{\sum_{i=1}^{m-1} \|\hat{\bx}_{t-i} - \hat{\bx}_{t-i+1} \|^2} } \|\hat{\epsilon}_t  \|.
\end{align*}
Using Taylor expansion, note we can also write $\hat{\epsilon}_t$ as 
\begin{align*}
    \hat{\epsilon}_t := & \ \begin{bmatrix}
        f(\hat{\bx}_{t-1}) - (f(\hat{\bx}_t) + \langle \nabla f(\hat{\bx}_t), \hat{\bx}_{t-1} - \hat{\bx}_t\rangle)  \\
        f(\hat{\bx}_{t-2}) - \left(f(\hat{\bx}_{t-1}) + \langle \nabla f(\hat{\bx}_t), \hat{\bx}_{t-2} - \hat{\bx}_{t-1}\rangle\right)  \\
        \vdots \\
        f(\hat{\bx}_{t-m+1}) - \left(f(\hat{\bx}_{t-m+2}) + \langle \nabla f(\hat{\bx}_t), \hat{\bx}_{t-m+1} - \hat{\bx}_{t-m+2}\rangle\right) \\
    \end{bmatrix} \\
            = & \ \begin{bmatrix}
        \langle \nabla f(\hat{\bx}_t), \hat{\bx}_{t-1} - \hat{\bx}_t \rangle + \sO(\frac{L}{2})\| \hat{\bx}_{t-1} - \hat{\bx}_t\|^2 - (\langle \nabla f(\hat{\bx}_t), \hat{\bx}_{t-1} - \hat{\bx}_t\rangle)  \\
        \langle \nabla f(\hat{\bx}_{t-1}), \hat{\bx}_{t-2} - \hat{\bx}_{t-1} \rangle + \sO(\frac{L}{2})\| \hat{\bx}_{t-2} - \hat{\bx}_{t-1}\|^2 - (\langle \nabla f(\hat{\bx}_t), \hat{\bx}_{t-2} - \hat{\bx}_{t-1}\rangle)   \\
        \vdots \\
 \langle \nabla f(\hat{\bx}_{t-m+2}), \hat{\bx}_{t-m+1} - \hat{\bx}_{t-m+2} \rangle + \sO(\frac{L}{2})\| \hat{\bx}_{t-m+1} - \hat{\bx}_{t-m+2}\|^2 - (\langle \nabla f(\hat{\bx}_t), \hat{\bx}_{t-m+1} - \hat{\bx}_{t-m+2}\rangle)  \\
    \end{bmatrix} \\
    = & \ \begin{bmatrix}
        \langle \nabla f(\hat{\bx}_t) - \nabla f(\hat{\bx}_t), \hat{\bx}_{t-1} - \hat{\bx}_t \rangle + \sO(\frac{L}{2})\| \hat{\bx}_{t-1} - \hat{\bx}_t\|^2  \\
        \langle \nabla f(\hat{\bx}_{t-1}) - \nabla f(\hat{\bx}_t) , \hat{\bx}_{t-2} - \hat{\bx}_{t-1} \rangle + \sO(\frac{L}{2})\| \hat{\bx}_{t-2} - \hat{\bx}_{t-1}\|^2    \\
        \vdots \\
 \langle \nabla f(\hat{\bx}_{t-m+2}) - \nabla f(\hat{\bx}_t), \hat{\bx}_{t-m+1} - \hat{\bx}_{t-m+2} \rangle + \sO(\frac{L}{2})\| \hat{\bx}_{t-m+1} - \hat{\bx}_{t-m+2}\|^2  \\
    \end{bmatrix} \\
    = 
    & \ \begin{bmatrix}
\sO(\frac{L}{2})\| \hat{\bx}_{t-1} - \hat{\bx}_t\|^2  \\
           \sO(L \|\hat{\bx}_{t-1}) - \nabla f(\hat{\bx}_t) \| \|\hat{\bx}_{t-2} - \hat{\bx}_{t-1} \|) + \sO(\frac{L}{2})\| \hat{\bx}_{t-2} - \hat{\bx}_{t-1}\|^2    \\
        \vdots \\
 \sO(L \|\hat{\bx}_{t-m+1}) - \nabla f(\hat{\bx}_t) \| \|\hat{\bx}_{t-m+1} - \hat{\bx}_{t-m+2} \|)  + \sO(\frac{L}{2})\| \hat{\bx}_{t-m+1} - \hat{\bx}_{t-m+2}\|^2  \\
    \end{bmatrix}.
\end{align*}
Thus it follows that 
\begin{align*}
    \|\hat{\epsilon}_t\| \leq \max \{\max_{i=1:m-1} \|\hat{\bx}_{t-i} - \hat{\bx}_t\|, \max_{i=1:m-1} \|\hat{\bx}_{t-i} - \hat{\bx}_{t-i+1}\|\} \sqrt{\sum_{i=1}^{m-1} \|\hat{\bx}_{t-i} - \hat{\bx}_{t-i+1} \|^2}
\end{align*}

It follows then that 
\begin{align*}
    \|\bm{g}_t - \nabla f(\hat{x}_t)\| \leq & \ \frac{1}{\sigma_{\min}(D_t)} \|\hat{\epsilon}_t \| \leq \frac{\sO(\sqrt{d}\sqrt{m})}{\sqrt{\sum_{i=1}^{m-1} \|\hat{\bx}_{t-i} - \hat{\bx}_{t-i+1} \|^2} } \|\hat{\epsilon}_t  \| \\
    \leq & \ \sO(\sqrt{d}\sqrt{m}) \max \{\max_{i=1:m-1} \| \hat{\bx}_{t-i} - \hat{\bx}_t\|, \max_{i=1:m-1} \|\hat{\bx}_{t-i} - \hat{\bx}_{t-i+1} \|\},
\end{align*}
i.e. the relation $C_{d,m} = \sO(\sqrt{d}\sqrt{m})$ holds.\footnote{We note that this is slightly different from what we stated in Assumption \ref{assumption:key_analysis_condition}, but this result is identical to that in Assumption \ref{assumption:key_analysis_condition} for the purposes of our proofs.}
\end{proof}

 \subsection{Proofs for Lemma \ref{lemma:xi_smaller_than_grad_per_step}, Theorem \ref{theorem:LReSZO-nonconvex-convergence}, and Theorem \ref{theorem:LReSZO-strong-convex-convergence}}
\label{app:Proofs_key_results}
    Throughout the analysis, without loss of generality, we denote $C_d:= C_{d,m}$ and assume that $C_d \geq 1$. 

We first prove the following useful helper result, which bounds the change in gradient norm across $m$ iterations whenever the gradient estimation bias can be bounded in terms of the true gradient norm.

\begin{lemma}
\label{lemma:gradient_norm_bounded_across_m_iters}
Suppose Assumption \ref{assumption:f_regularity} holds. Suppose $m > d$. 
Consider any $t\geq m$. Suppose that 
\begin{align}
\label{eq:1_step_bdd}
    \|\bx_{j+1} - \bx_j\| \leq 2\eta \|\nabla f(\bx_j)\|  
\end{align}
for all $j < t$.

Then, by picking $\eta = \sO\left(\frac{1}{m L} \right)$, we can ensure that 
\begin{align*}
    \|\nabla f(\bx_{t-j}) \| \leq 4 \|\nabla f(\bx_{t})\|
\end{align*}
for any $j \in [m]$.
\end{lemma}
\begin{proof}
Consider any $t \geq m$, and consider any $0 \leq j < t$. Since by assumption
\begin{align*}
    \|\bx_{j+1} - \bx_j\| \leq 2\eta \|\nabla f(\bx_j)\|,
\end{align*}
we have that
\begin{align}
    \|\nabla f(\bx_{j+1}) - \nabla f(\bx_j)\| \leq \ L \|\bx_{j+1} - \bx_j \| \leq \ 2\eta L \|\nabla f(\bx_j)\|, \label{eq:grad_f_1_step_bdd_by_assumption}
\end{align}
where the first line utilizes $L$-smoothness and the second line utilizes our assumption in \eqref{eq:1_step_bdd}. Continuing, we observe that
\begin{align*}
    \|\nabla f(\bx_j)\| \leq & \ \|\nabla f(\bx_{j+1}) \| + \|\nabla f(\bx_j) - \nabla f(\bx_{j+1}) \| \\
    \leq & \|\nabla f(\bx_{j+1}) \| + 2\eta L \|\nabla f(\bx_j) \|, 
\end{align*}
where the first inequality uses triangle inequality and the second inequality follows from \eqref{eq:grad_f_1_step_bdd_by_assumption}. Rearranging, we find that
\begin{align*}
    & \ \left(1 - 2 \eta L\right) \|\nabla f(\bx_j) \| \leq \| \nabla f(\bx_{j+1})\| \\
    \implies & \  \|\nabla f(\bx_j) \| \leq \frac{1}{1 - 2 \eta L} \| \nabla f(\bx_{j+1})\|.
\end{align*}
Suppose we pick $\eta = \frac{1}{2 m L}$. Then, continuing, we find that
\begin{align*}
    \|\nabla f(\bx_j) \| \leq & \ \frac{1}{1 - 2 \eta L} \| \nabla f(\bx_{j+1})\| \\
    \leq & \ \frac{1}{1- \frac{1}{m}} \|\nabla f(\bx_{j+1})\| \\
    = & \ \frac{m}{m-1} \|\nabla f(\bx_{j+1})\|.
\end{align*}
It follows then that for any $t \geq m$, for any $j \in [m]$, we have
\begin{align*}
    \|\nabla f(\bx_{t-j}) \| \leq & \ \left(\frac{m}{m-1}\right)^j \|\nabla f(\bx_{t})\| \\
    \leq & \ \left(\frac{m}{m-1}\right)^m \|\nabla f(\bx_{t})\| \\
    \leq & \ 4 \|\nabla f(\bx_{t})\|,
\end{align*}
where the final inequality follows from the fact that the sequence $(\frac{n}{n-1})^n \leq 4$ for all $n \geq 2$. 
\end{proof}

Next, we restate and prove Lemma \ref{lemma:xi_smaller_than_grad_per_step}, which shows that for any $0 < c < 1$, we can inductively guarantee that the gradient estimation bias $\|\xi_t \|$ always satisfies $\|\xi_t\| \leq c \|\nabla f(\bx_t)\|$ for all $t \geq m$. 

\lemmaXiSmallerThanGrad*
\begin{proof}
We consider an inductive proof. Our inductive claim is as follows. Pick any $t \geq m$. suppose that the inductive assumption
\begin{align*}
        \|\bx_{j+1} - \bx_j\| \leq 2\eta \| \nabla f(\bx_j)\|
\end{align*}
holds for all $j < t$. Then, for any $0 < c <1$, by picking $\eta = C\frac{1}{m C_d L}$ for some sufficiently small absolute constant $C > 0$, we will show that the inductive claim implies that
\begin{align*}
    & \ \|\xi_t\| \leq c \|\nabla f(\bx_t)\| \\
    \implies & \ \|\bx_{t+1} - \bx_t\| = \| \eta \bm{g}_t\| \leq \eta(\|\nabla f(\bx_t)\| + \| \xi_t\|)  \leq 2\eta \| \nabla f(\bx_t)\|,
\end{align*}
which inductively guarantees that $\|\xi_t \| \leq c \| \nabla f(\bx_t)\|$ for all $t'\geq t$.

We first observe that our inductive assumption holds for the base case $t = m$, since by our assumption, 
\begin{align*}
        \|\bx_{j+1} - \bx_j\| \leq 2\eta \| \nabla f(\bx_j)\|
\end{align*}
holds for all $j < m$.

Let $\eta_t$ denote the step-size at time $t$, which is $\tilde{\eta}$ for any $t < m$ and $\eta$ after that. Observe that for any $t \geq m$,\footnote{For notational simplicity, our analysis below assumes that the first $m$ iterations also utilizes the adaptive smoothing scheme. However, this is by no means necessary and we note that with some more tedious derivations, we can easily generalize the derivation to the case with no adaptive smoothing in the first $m$ iterations and achieve the same bound.}
\small
\begin{align}
        &  \ \|\xi_t\| \leq  \ \|\bm{g}_t - \nabla f(\hat{\bx}_t) \| + \|\nabla f(\hat{\bx}_t) - \nabla f(\bx_t )\| \nonumber \\
        \labelrel{\leq}{eq:xi_bdd_L_smooth_second} & \ \|\bm{g}_t - \nabla f(\hat{\bx}_t) \| + L \|\hat{\bx}_t - \bx_t\| \nonumber \\
  \labelrel{\leq}{eq:xi_bdd_xhat_defn_second} & \ \|\bm{g}_t - \nabla f(\hat{\bx}_t) \| + \eta_{t-1} L \|\bm{g}_{t-1} \| \nonumber \\
        \labelrel{\leq}{eq:xi_bdd_assumption_2_second} & \ C_d \frac{L}{2} \max_{i=1:m-1}\left\|\hat{\bx}_{t-i} - \hat{\bx}_t\right\| + \eta_{t-1} L \|\bm{g}_{t-1} \| \nonumber \\
        = & \ \frac{C_d L}{2} \max_{i=1:m-1}\left\|\sum_{j=1}^i (\hat{\bx}_{t-j} - \hat{\bx}_{t-j+1})\right\| + \eta_{t-1} L \|\bm{g}_{t-1} \| \nonumber \\
        = & \ \frac{C_d L}{2} \max_{i=1:m-1}\left\|\sum_{j=1}^i (\bx_{t-j} - \bx_{t-j+1} + \delta_{t-j} \bu_{t-j} - \delta_{t-j+1} \bu_{t-j+1})\right\| + \eta_{t-1} L \|\bm{g}_{t-1} \| \nonumber \\
        = & \ \frac{C_d L}{2} \max_{i=1:m-1}\left\|\sum_{j=1}^i (\bx_{t-j} - \bx_{t-j+1} + \eta_{t-j-1} \|\bm{g}_{t-j-1}\| \bu_{t-j} - \eta_{t-j} \|\bm{g}_{t-j}\|\bu_{t-j+1})\right\| \nonumber  \\
        & \ \quad \quad + \eta_{t-1} L \|\bm{g}_{t-1} \| \nonumber \\
         \labelrel{\leq}{eq:xi_bdd_u_unit_second} &  \ \frac{C_d L}{2}\max_{i=1:m-1} \left( \sum_{j=1}^i \left(\eta_{t-j}\left\| \bm{g}_{t-j}\right\|  + \eta_{t-j-1} \| \bm{g}_{t-j-1}\| + \eta_{t-j} \|\bm{g}_{t-j}\|\right)\right) + \eta_{t-1} L \|\bm{g}_{t-1} \| \nonumber \\
        \leq & \ \frac{C_d L}{2}\max_{i=1:m-1} \left( \sum_{j=1}^i \left(2\eta_{t-j} \|\bm{g}_{t-j}\|  + \eta_{t-j-1} \| \bm{g}_{t-j-1}\|\right)\right) + \eta_{t-1} L \|\bm{g}_{t-1} \| \nonumber \\
    \labelrel{\leq}{eq:xi_bdd_Cd_larger_than_1_second} & \  \frac{C_d L}{2}\left( \sum_{j=1}^{m-1} \left(4\eta_{t-j} \|\bm{g}_{t-j}\|  + \eta_{t-j-1} \| \bm{g}_{t-j-1}\| \right)\right) \nonumber \\
    \leq & \ \frac{C_d L}{2}\sum_{j=1}^{m} 4\eta_{t-j} \|\bm{g}_{t-j}\| \nonumber \\
    \labelrel{=}{eq:x_1_step_eta_g_defn} & \ \frac{C_d L}{2}\sum_{j=1}^{m} 4\| \bx_{t-j+1} - \bx_{t-j}\| \nonumber\\
    \labelrel{\leq}{eq:inductive_assumption} & \ 2C_d L \sum_{j=1}^{m} 2 \eta \|\nabla f(\bx_{t-j})\|, \nonumber
\end{align}
\normalsize
where    \eqref{eq:xi_bdd_L_smooth_second} holds by $L$-smoothness of $f$, \eqref{eq:xi_bdd_xhat_defn_second} holds by definition of $\hat{\bx}_t := \hat{\bx}_t + \eta_{t-1} \|\bm{g}_{t-1}\| \bu_t$ due to the adaptive selection of the smoothing radius $\delta_t := \eta \|\bm{g}_{t-1}\|$ and the fact that each $\bu_t$ has unit norm, \eqref{eq:xi_bdd_assumption_2_second} utilizes Assumption \ref{assumption:key_analysis_condition}, \eqref{eq:xi_bdd_u_unit_second} again uses the fact that each $\bu_t$ has a unit norm, \eqref{eq:xi_bdd_Cd_larger_than_1_second} holds by our assumption at the start of this subsection that $C_d \geq 1$, \eqref{eq:x_1_step_eta_g_defn} holds since $\bx_{t-j+1} = \bx_{t-j} + \eta_{t-j} \bm{g}_{t-j}$, and \eqref{eq:inductive_assumption} follows by our inductive assumption.

Now, utilizing our inductive assumption, we may also use Lemma \ref{lemma:gradient_norm_bounded_across_m_iters} to find that for any $j \in [m]$, $\|\nabla f(\bx_{t-j})\| \leq 4 \|\nabla f(\bx_t) \|$. Thus, continuing from above, we find that
\begin{align*}
    \|\xi_t\| \leq & \  2 C_d L \sum_{j=1}^{m} 2 \eta \|\nabla f(\bx_{t-j})\| \\
    \leq & \ 4 C_d L \sum_{j=1}^{m} 4 \eta \|\nabla f(\bx_{t})\| \\
    \leq & \ 16 C_d L \eta \sum_{j=1}^{m} \|\nabla f(\bx_{t})\|.
\end{align*}
Now, by picking $\eta  = \frac{C}{m C_d L}$ for a sufficiently small $C > 0$, we may conclude that 
\begin{align*}
    & \ \|\xi_t\| \leq c \|\nabla f(\bx_t) \|, \\
    \implies & \ \|\bx_{t+1} - \bx_t\| = \| \eta \bm{g}_t\| \leq \eta(\|\nabla f(\bx_t)\| + \| \xi_t\|)  \leq 2\eta \| \nabla f(\bx_t)\|.
\end{align*}
Thus, by induction, since the base case of $t = m$ holds, we can conclude that for all $t \geq m$, we have
\begin{align*}
    & \ \|\xi_t\| \leq c \|\nabla f(\bx_t) \|, 
\end{align*}
as desired.
\end{proof}

We next restate and prove Lemma \ref{lemma:fn_improvement_per_step}.

\fnImprovePerStep*

\begin{proof}
    Consider any positive integer $ t \geq m$. 
    Fix any $0 < c < 1$. By our assumptions, applying Lemma \ref{lemma:xi_smaller_than_grad_per_step}, we find that there exists a sufficiently small absolute constant $C > 0$ such that when $\eta = \frac{C}{m C_d L}$, we have
    \begin{align*}
        \|\xi_t\| \leq c \| \nabla f(\bx_t)\|
    \end{align*}
    for any $t \geq m$.

    Next, observe that
    \begin{align}
    & \ f(\bx_{t+1}) - f(\bx_t) \nonumber \\
    \labelrel{\leq}{eq:prop_1_L_smooth} & \  -\eta \langle \nabla f(\bx_t),  \bm{g}_t \rangle + \frac{L}{2} \|\eta \bm{g}_t \|^2  \nonumber \\
    = & \ -\eta \langle \nabla f(\bx_t),  \nabla f(\bx_t) \rangle - \eta \langle \nabla f(\bx_t),  \bm{g}_t - \nabla f(\bx_t) \rangle \nonumber \\
    & \quad \quad + \frac{L}{2} \|\eta (\bm{g}_t - \nabla f(\bx_t) + \nabla f(\bx_t)) \|^2  \nonumber \\
    \labelrel{\leq}{eq:prop_1_young} & \ -\eta \|\nabla f(\bx_t) \|^2 + \frac{\eta}{2} \left(\|\nabla f(\bx_t)\|^2 + \|\xi_t \|^2  \right) \nonumber \\
    & \quad \quad + \eta^2 L \left(\|\nabla f(\bx_t)\|^2 + \|\xi_t\|^2 \right) \nonumber  \\
    \leq & \ -\frac{\eta}{4} \|\nabla f(\bx_t) \|^2 + \eta \|\xi_t \|^2, \label{eq:f_decrease_1_step_L_smooth}
\end{align}
where \eqref{eq:prop_1_L_smooth} follows by $L$-smoothness of $f$, \eqref{eq:prop_1_young} follows by application of Young's equality as well as the notation $\xi_t := \bm{g}_t - \nabla f(\bx_t)$, and the final inequality follows by choosing $\eta \leq \frac{1}{4L}$. 

     Thus, continuing from \eqref{eq:f_decrease_1_step_L_smooth}, by picking $c$ to be smaller than $\frac{1}{\sqrt{8}}$, such that
    \begin{align*}
        \|\xi_t\| \leq \frac{1}{\sqrt{8}}\| \nabla f(\bx_t)\|,
    \end{align*}
    we find that 
    \begin{align*}
            f(\bx_{t+1}) - f(\bx_t) 
    \leq & \ -\frac{\eta}{4} \|\nabla f(\bx_t) \|^2 + \eta \|\xi_t \|^2. \\
    \leq & \ -\frac{\eta}{4} \|\nabla f(\bx_t) \|^2 + \frac{\eta}{8}\|\nabla f(\bx_t) \|^2 \\
    \leq & \ - \frac{\eta}{8} \|\nabla f(\bx_t) \|^2,
    \end{align*} 
    as desired.

\end{proof}

We next restate and prove Theorem \ref{theorem:LReSZO-nonconvex-convergence}, which provides a convergence result for smooth nonconvex functions.

\thmNonconvex*

\begin{proof}
   By Lemma \ref{lemma:fn_improvement_per_step}, we have that
\begin{align*}
            f(\bx_{t+1}) - f(\bx_t) \leq - \frac{\eta}{8} \|\nabla f(\bx_t) \|^2. 
\end{align*}
Summing this from $t = m$ to $t = T$ and plugging in our choice of $\eta = \sO\left(\frac{1}{m C_d L}\right)$, we observe that
\begin{align*}
    & \ f(\bx_T) - f(\bx_m) \leq - \sO\left(\frac{1}{m C_dL }\sum_{t=m}^{T-1} \|\nabla f(\bx_t) \|^2\right) \\
    \implies & \ \sO\left(\frac{1}{m C_dL }\sum_{t=m}^{T-1} \|\nabla f(\bx_t)\|^2\right) \leq f(\bx_m) - f(\bx_T) \\
    \implies & \ \frac{1}{T}\sum_{t=m}^{T-1} \|\nabla f(\bx_t)\|^2 \leq \sO \left(\frac{m C_d L \left(f(\bx_m) - f^*\right)}{T}\right),
\end{align*}
where the final inequality uses the simple fact that $f(\bx_t) \geq f^*$.
\end{proof}

We next restate and prove Theorem \ref{theorem:LReSZO-strong-convex-convergence}, which shows a linear convergence rate for smooth and strongly convex functions.

\thmStrongConvex*

\begin{proof}
    Consider any positive integer $t \geq m$. 
By Lemma \ref{lemma:fn_improvement_per_step}, we have that
    \begin{align*}
            f(\bx_{t+1}) - f(\bx_t)  \leq & \ - \frac{\eta}{8} \|\nabla f(\bx_t) \|^2 \\
    \leq & \ - \frac{\eta\mu}{4} (f(\bx_t) - f^*),\\
    \Longrightarrow \  f(\bx_{t+1}) - f^* &\leq \left(1-\frac{\eta\mu}{4}\right) (f(\bx_t) - f^*), 
    \end{align*} 
    where the fourth inequality follows by $\mu$-strong convexity which ensures that $\| \nabla f(\bx_t)\|^2 \geq 2\mu(f(\bx_t) - f^*)$. Our desired result then follows. 
\end{proof}

\newpage 
\section{Empirical Study of \texorpdfstring{$C_d$}{Cd} in Assumption \ref{assumption:key_analysis_condition}}
We assume throughout this section that $m = \sO(d)$, so we may denote $C_d := C_{d,m}$.

We recall from Assumption \ref{assumption:key_analysis_condition} that the $C_d$ term is defined as a positive ratio term such that 
\begin{align*}
    \|\bm{g}_t - \nabla f(\hat{\bx}_t)\| \leq C_d \frac{L}{2} \max_{i=1:m-1}\|\hat{\bx}_{t-i} - \hat{\bx}_t\|.
\end{align*}
In our empirical calculations below, for simplicity, we calculate this ratio $C_d$ at the $t$-th iteration of the algorithm as 
\begin{align*}
C_d(t) := \frac{\|\bm{g}_t - \nabla f(\hat{\bx}_t) \|}{\frac{L}{2} \|\hat{\bx}_{t-(m-1)} - \hat{\bx}_t\|} \geq \frac{\|\bm{g}_t - \nabla f(\hat{\bx}_t) \|}{\frac{L}{2} \max_{i=1:m-1}\|\hat{\bx}_{t-i} - \hat{\bx}_t\|},
\end{align*}
and use our analytical knowledge of the ridge and logistic regression problems in \eqref{eq:ridge} and \eqref{eq:logis} to compute the smoothness constant $L$.

\subsection{Ridge Regression}

\begin{figure}[ht]
  \centering
  \begin{subfigure}[b]{0.32\textwidth}
    \includegraphics[width=\linewidth]{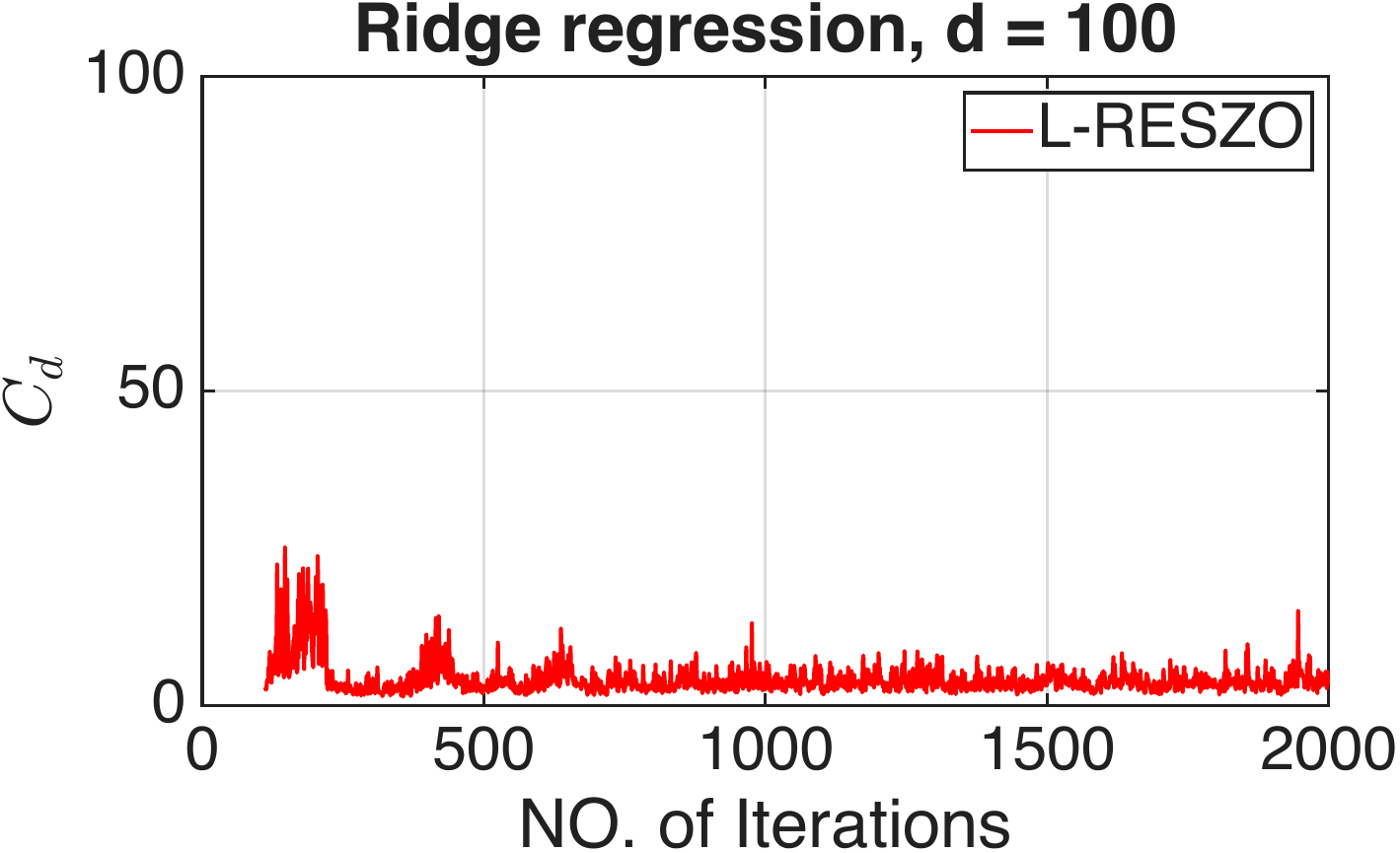}
    \caption{$d=100$}
    \label{fig:cd100_ridge}
  \end{subfigure}
  \hfill
  \begin{subfigure}[b]{0.32\textwidth}
    \includegraphics[width=\linewidth]{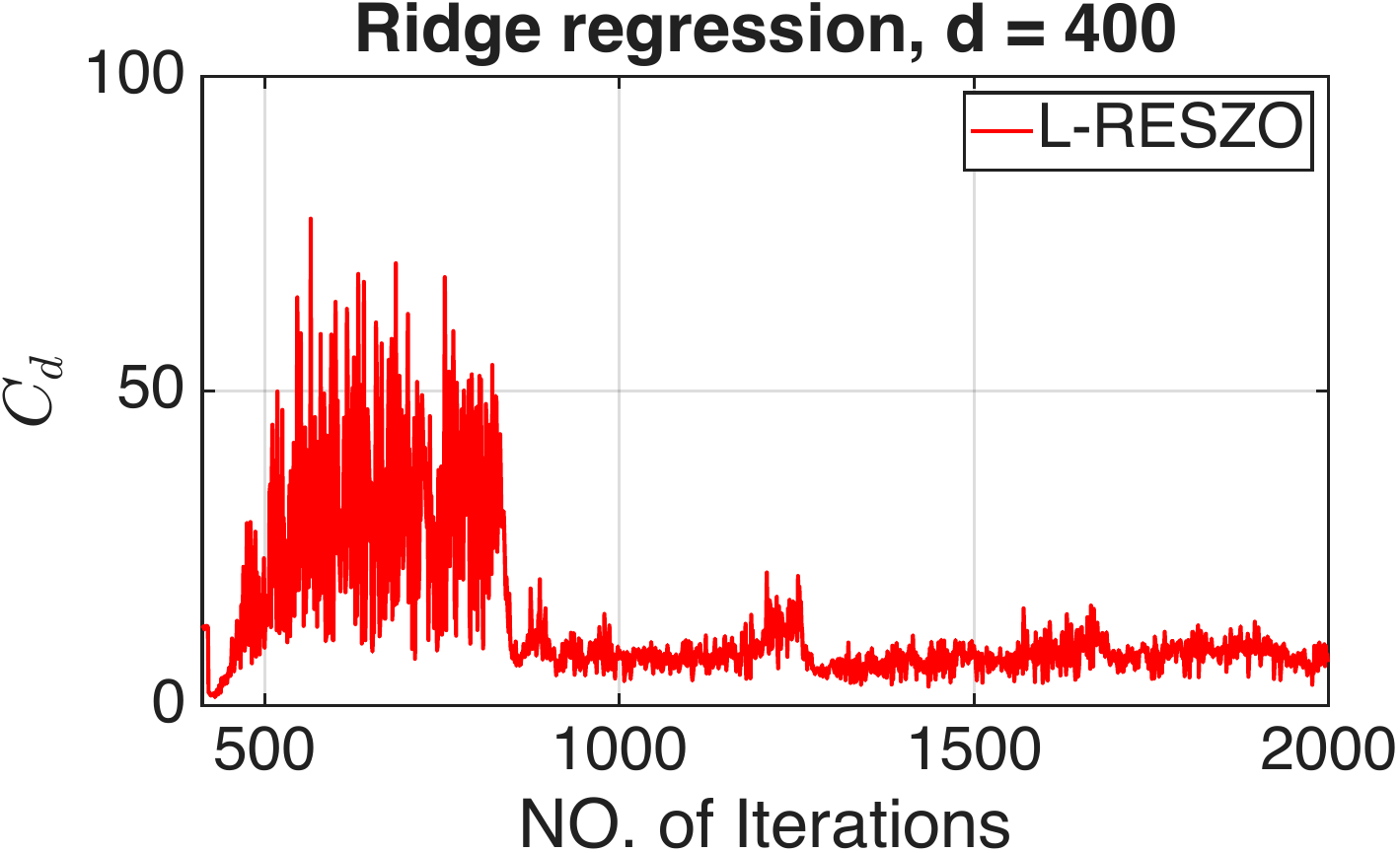}
    \caption{$d=400$}
    \label{fig:cd400_ridge}
  \end{subfigure}
  \hfill
  \begin{subfigure}[b]{0.32\textwidth}
    \includegraphics[width=\linewidth]{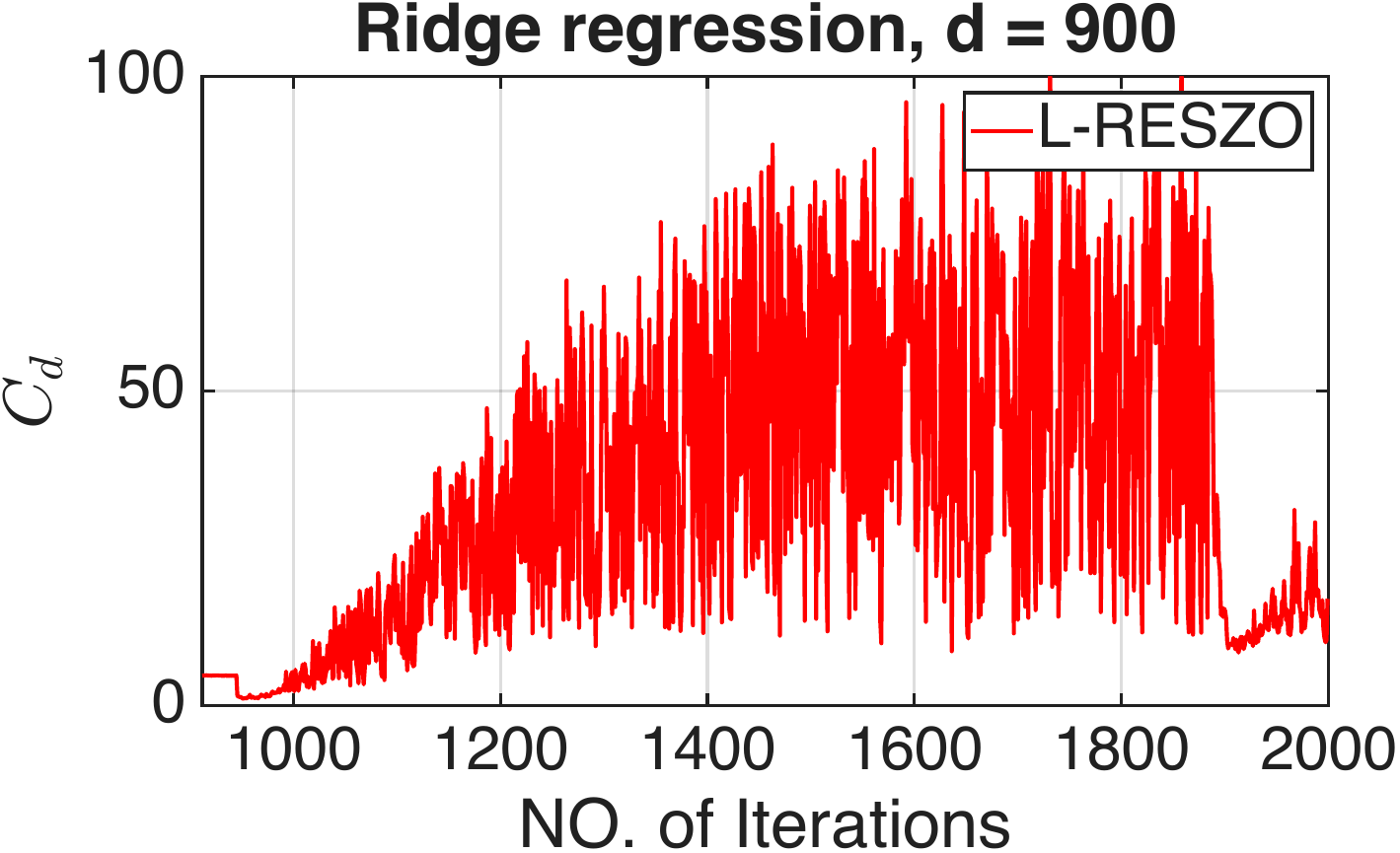}
    \caption{$d=900$}
    \label{fig:cd900_ridge}
  \end{subfigure}

  \vspace{0.5em}
  \caption{Empirical distribution of the ratio $C_d$ for dimensions
           $d=100,400,900$ for ridge regression problem in \eqref{eq:ridge}.  These plots are
           the sources from which the max/99-percentile/mean in
           Table~\ref{tab:ridge_Cd} are computed.}
  \label{fig:cd_plots_ridge}
\end{figure}

\begin{figure}[ht]
  \centering
  \begin{subfigure}[b]{0.32\textwidth}
    \includegraphics[width=\linewidth]{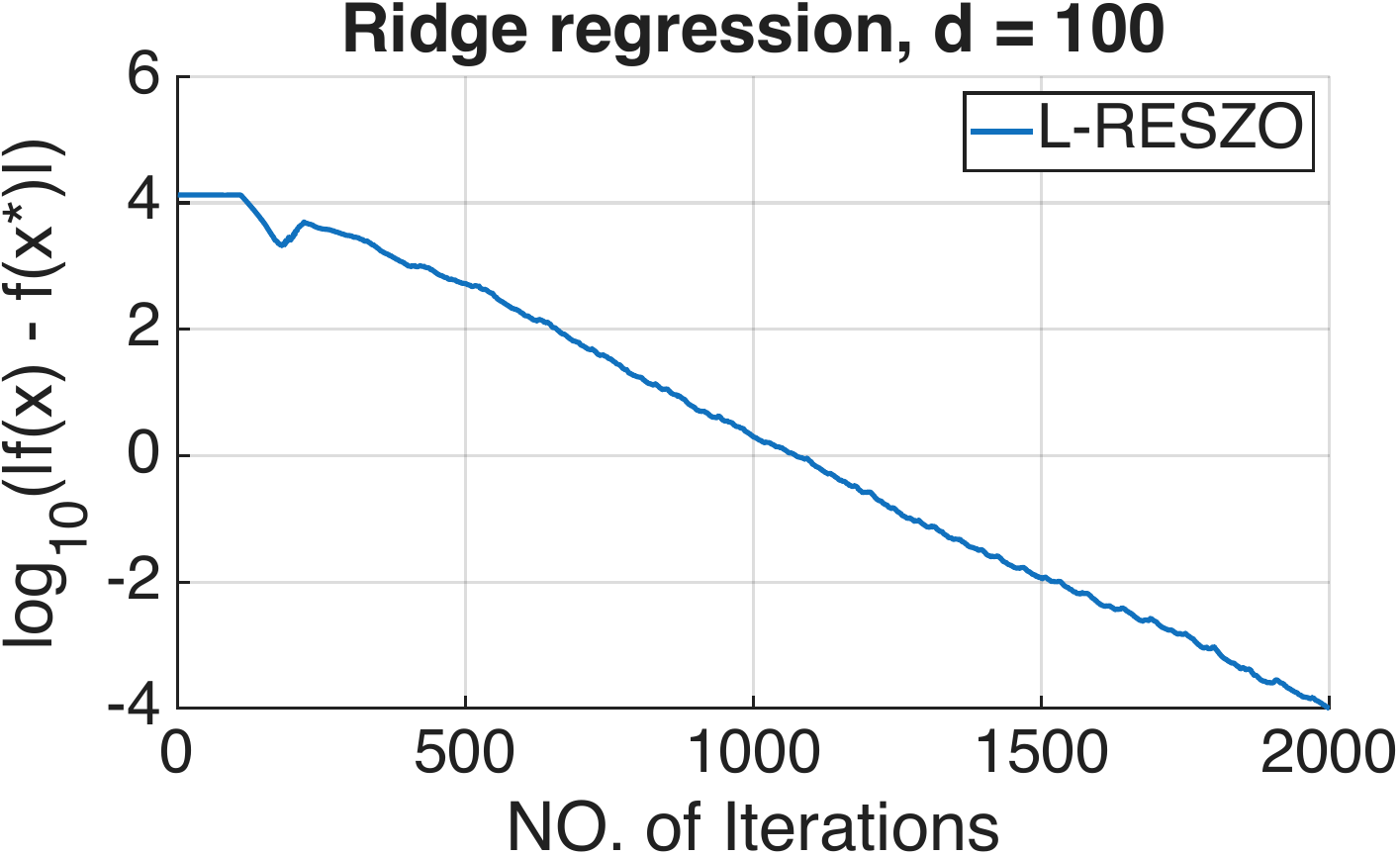}
    \caption{\(d=100\)}
    \label{fig:conv_ridge_100}
  \end{subfigure}
  \hfill
  \begin{subfigure}[b]{0.32\textwidth}
    \includegraphics[width=\linewidth]{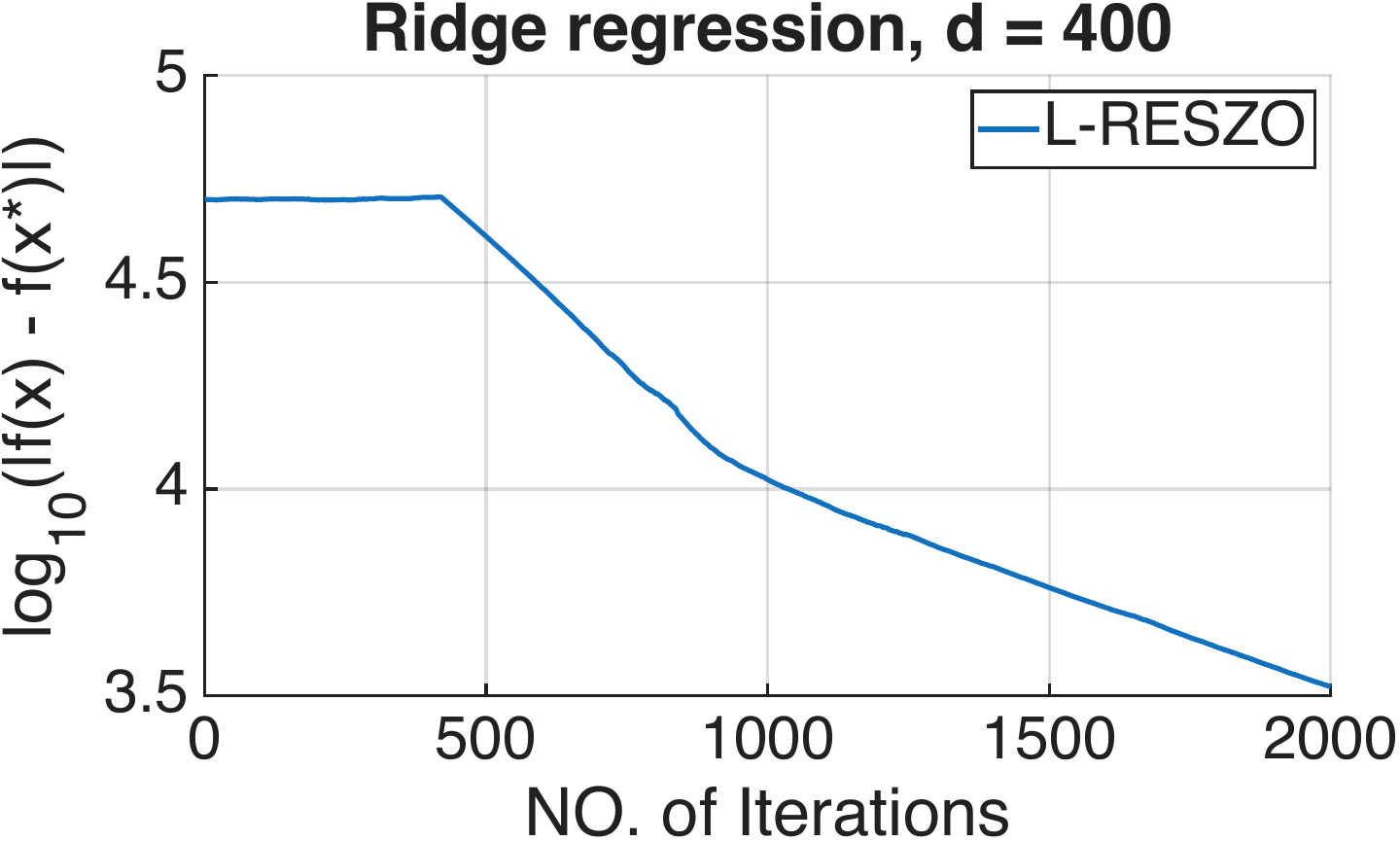}
    \caption{\(d=400\)}
    \label{fig:conv_ridge_400}
  \end{subfigure}
  \hfill
  \begin{subfigure}[b]{0.32\textwidth}
    \includegraphics[width=\linewidth]{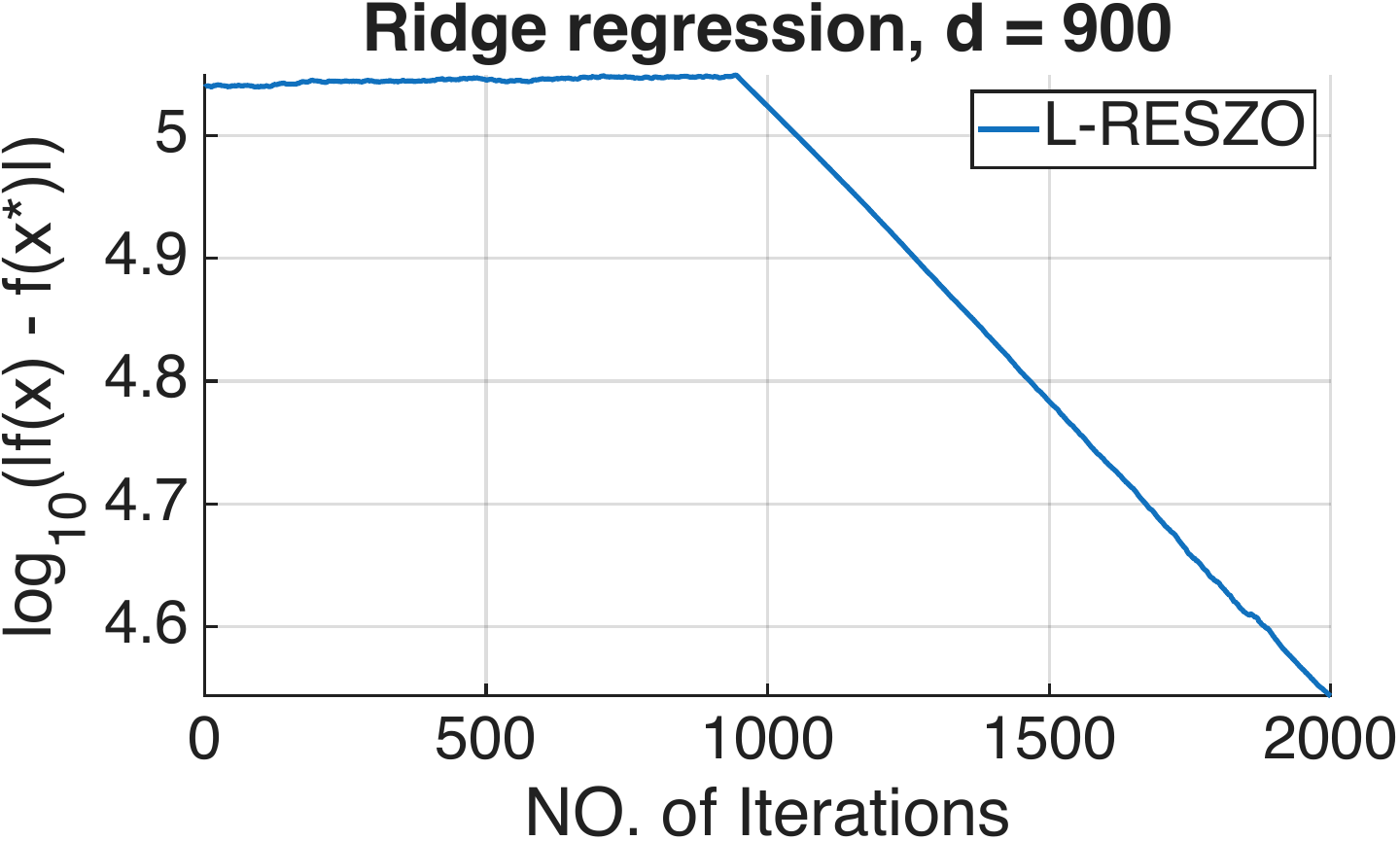}
    \caption{\(d=900\)}
    \label{fig:conv_ridge_900}
  \end{subfigure}

  \vspace{0.5em}
  \caption{Convergence of L-ReSZO on the ridge‐regression problem in \eqref{eq:ridge} for dimensions
           \(d=100,400,900\). The $C_d$ ratios shown in Figure \ref{fig:cd_plots_ridge} are based on the algorithmic trajectories here.}
  \label{fig:conv_ridge_all}
\end{figure}

We experimented on the ridge regression problem in \eqref{eq:ridge} for three values for the dimension $d$: $d = 100, 400, 900$. We computed the $C_d$ ratio statistics for a trial of the algorithm for these three dimensions, and display this information in Table \ref{tab:ridge_Cd}, which is based on the data from Figure \ref{fig:cd_plots_ridge}; the corresponding convergence behavior is shown in Figure \ref{fig:conv_ridge_all}. We observe that the maximum $C_d$ scales as around $3\sqrt{d}$, i.e., $\sO(\sqrt{d})$.

\begin{table}[ht]
  \captionsetup{width=\textwidth}
  \centering
  \caption{Statistics of the $C_d$ ratio across different dimensions $d$ for the ridge regression problem  \eqref{eq:ridge}.}
  \label{tab:ridge_Cd}
  \begin{tabular*}{\textwidth}{@{\extracolsep{\fill}} c c c c @{}}
    \toprule[1.2pt]
    \textbf{Dimension $d$} & \textbf{Max} & \textbf{99th Percentile} & \textbf{Mean} \\
    \midrule
    $100$ & $25.1522$  & $14.9690$  & $3.9955$  \\
    $400$ & $77.3924$  & $53.6369$  & $10.4881$ \\
    $900$ & $104.1276$ & $85.5617$  & $18.3887$ \\
    \bottomrule[1.2pt]
  \end{tabular*}
\end{table}

\subsection{Logistic Regression}

\begin{figure}[H]
  \centering
  \begin{subfigure}[b]{0.32\textwidth}
    \includegraphics[width=\linewidth]{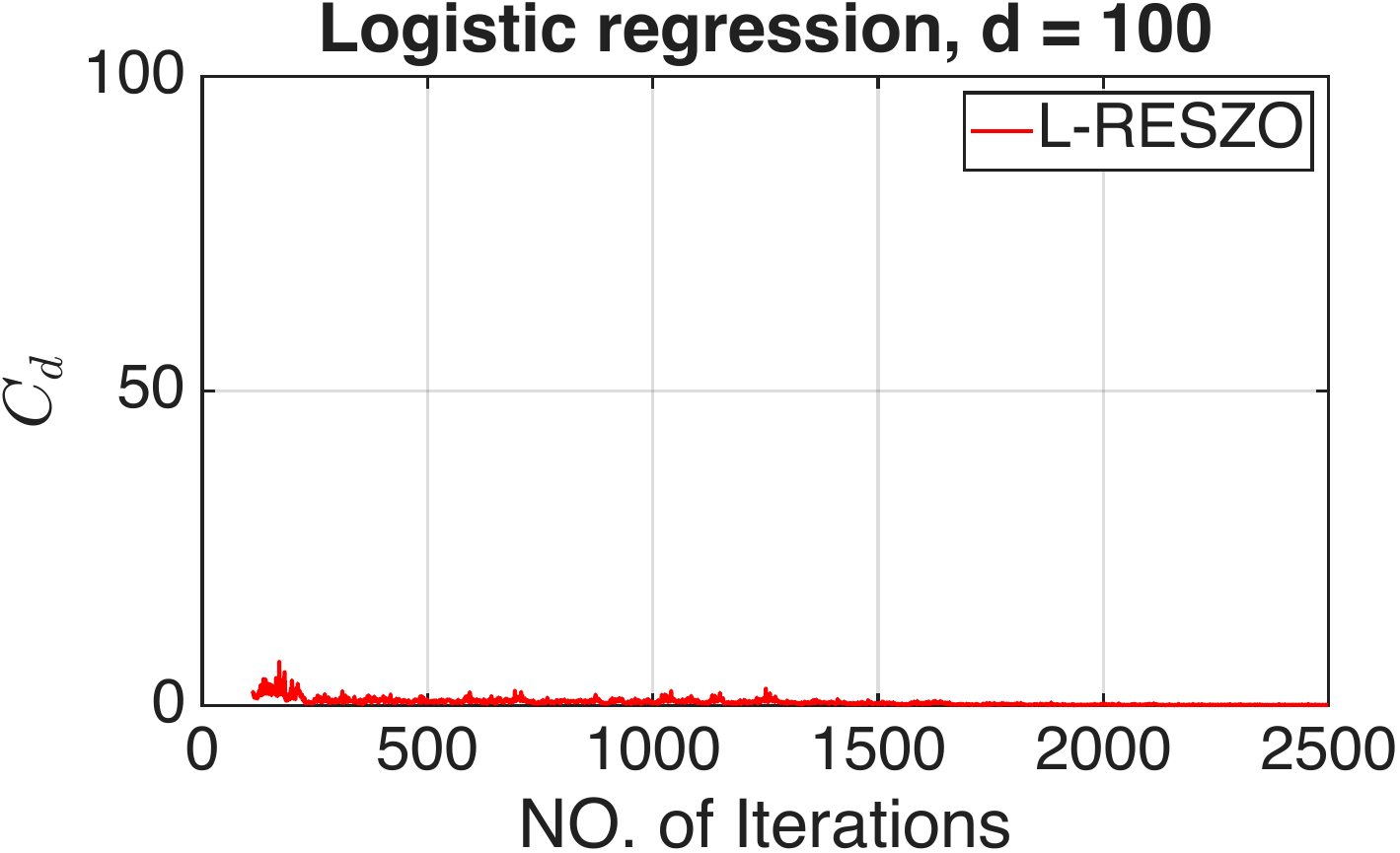}
    \caption{$d=100$}
    \label{fig:cd100}
  \end{subfigure}
  \hfill
  \begin{subfigure}[b]{0.32\textwidth}
    \includegraphics[width=\linewidth]{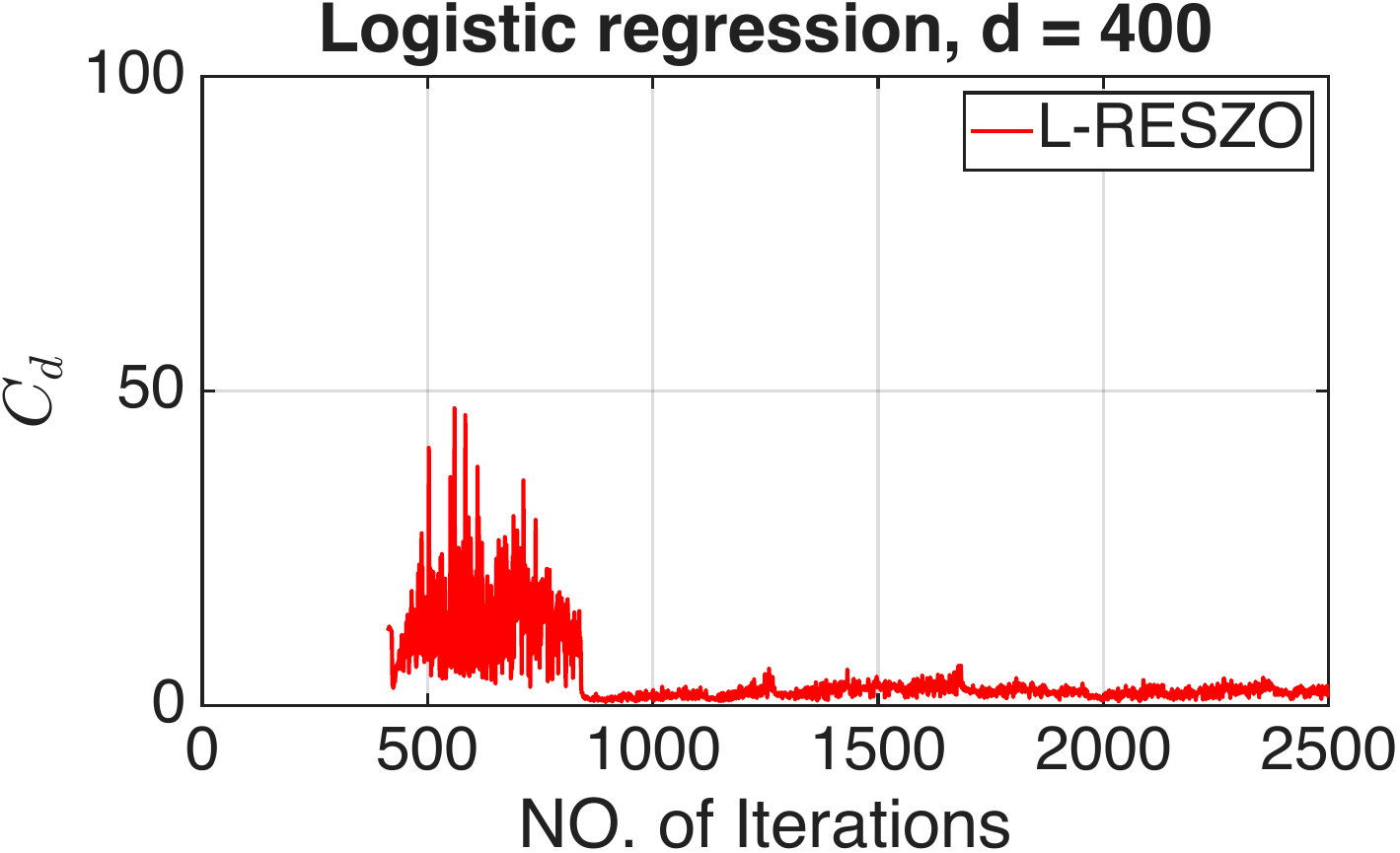}
    \caption{$d=400$}
    \label{fig:cd400}
  \end{subfigure}
  \hfill
  \begin{subfigure}[b]{0.32\textwidth}
    \includegraphics[width=\linewidth]{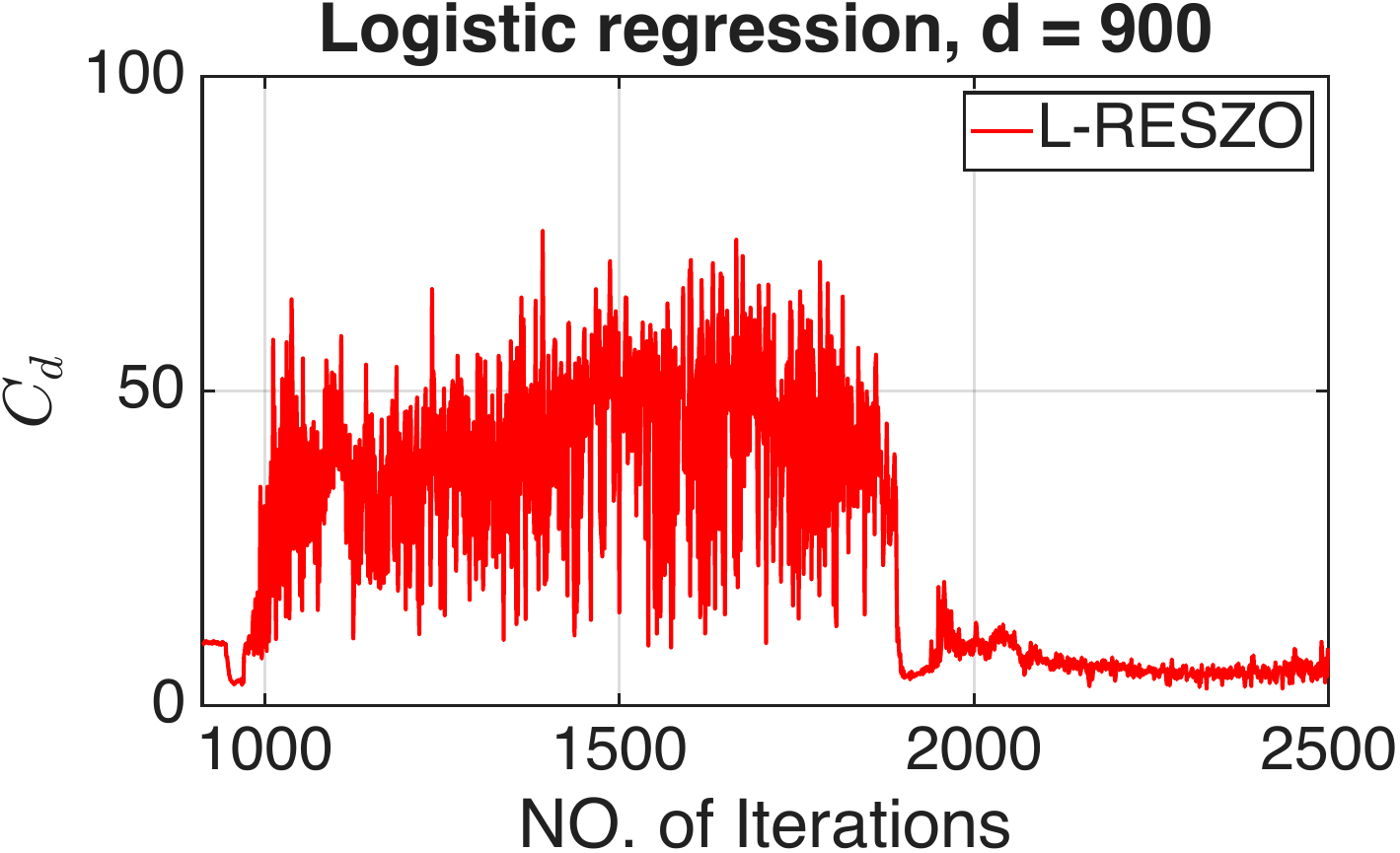}
    \caption{$d=900$}
    \label{fig:cd900}
  \end{subfigure}

  \vspace{0.5em}
  \caption{Empirical distribution of the ratio $C_d$ for dimensions
           $d=100,400,900$ for logistic regression problem in \eqref{eq:logis}.  These plots are
           the sources from which the max/99-percentile/mean in
           Table~\ref{tab:logistic_Cd} were computed.}
  \label{fig:cd_plots_logistics}
\end{figure}

\begin{figure}[H]
  \centering
  \begin{subfigure}[b]{0.32\textwidth}
    \includegraphics[width=\linewidth]{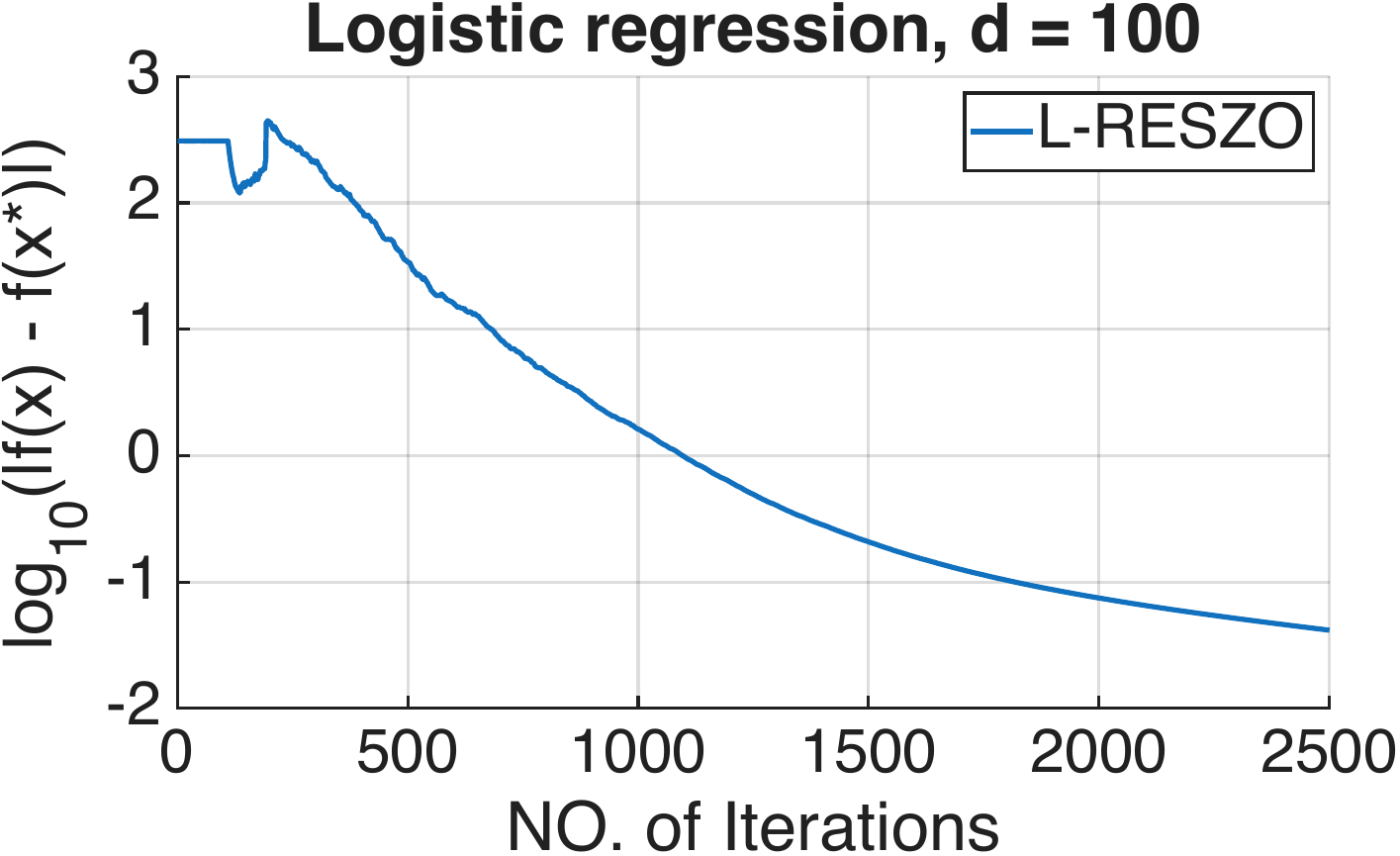}
    \caption{\(d=100\)}
    \label{fig:conv_logis_100}
  \end{subfigure}
  \hfill
  \begin{subfigure}[b]{0.32\textwidth}
    \includegraphics[width=\linewidth]{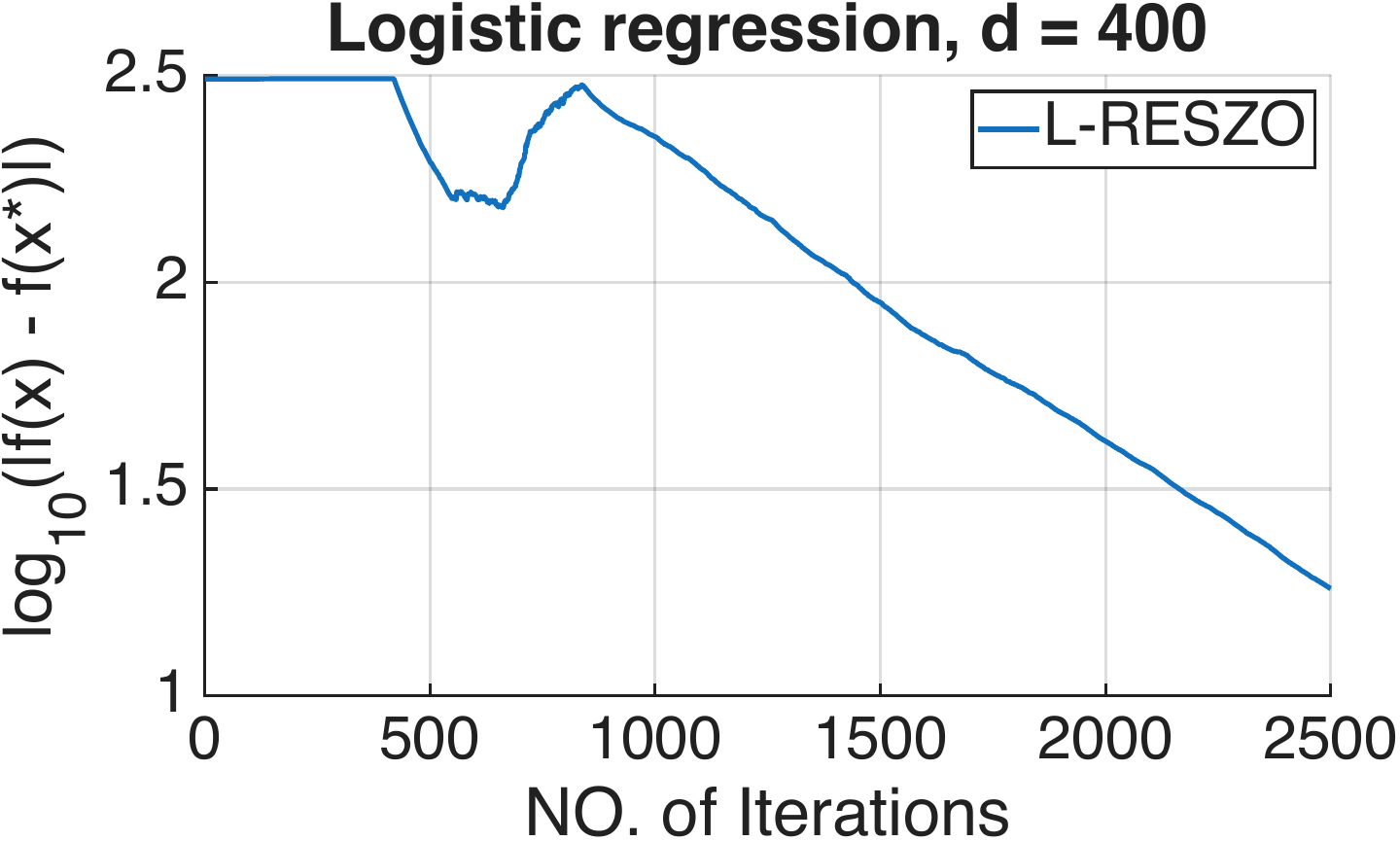}
    \caption{\(d=400\)}
    \label{fig:conv_logis_400}
  \end{subfigure}
  \hfill
  \begin{subfigure}[b]{0.32\textwidth}
    \includegraphics[width=\linewidth]{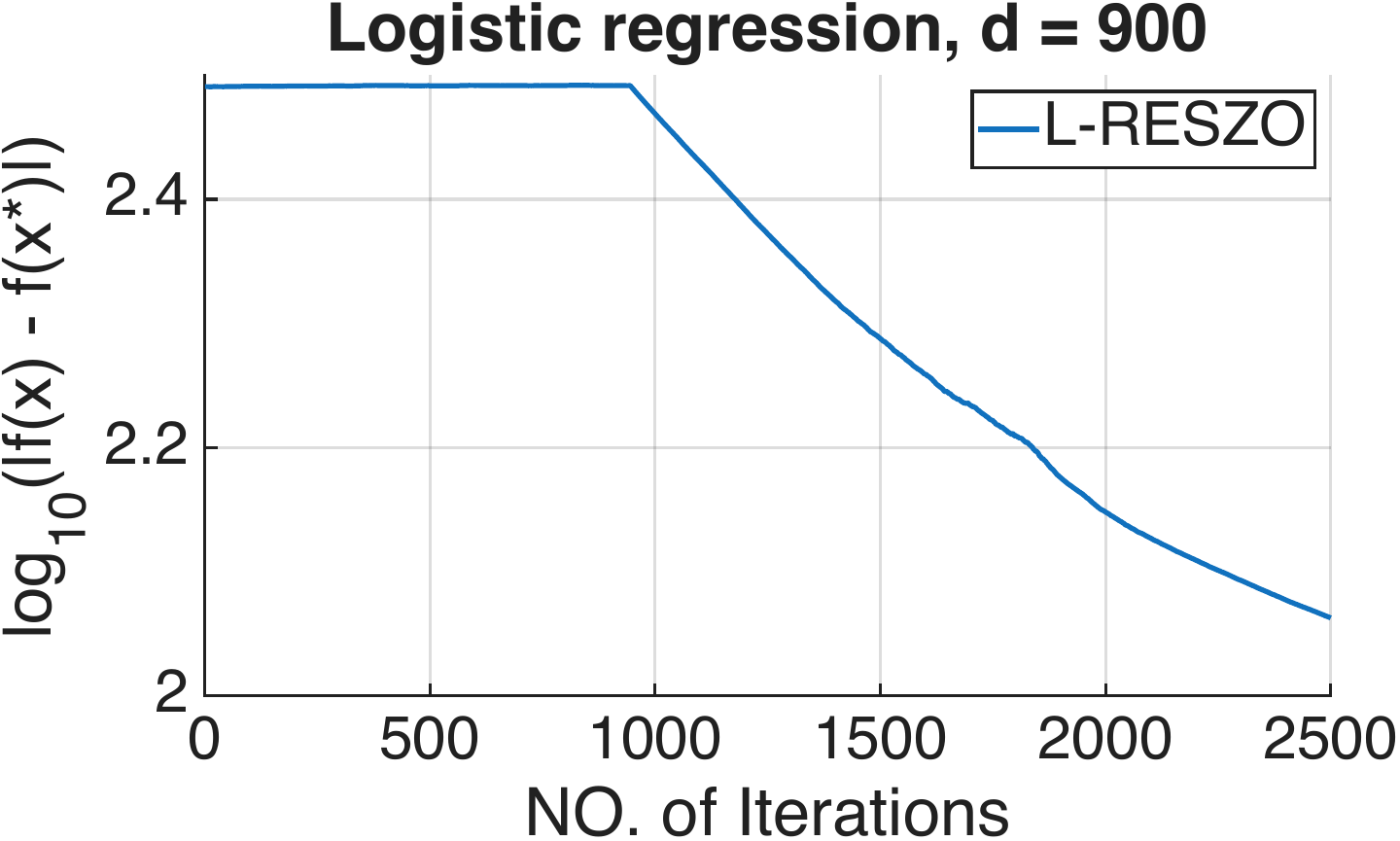}
    \caption{\(d=900\)}
    \label{fig:conv_logis_900}
  \end{subfigure}
  \vspace{0.5em}
  \caption{Convergence of L-ReSZO on the logistic regression problem in \eqref{eq:logis} for dimensions
           \(d=100,400,900\). The $C_d$ ratios shown in Figure \ref{fig:cd_plots_logistics} are based on the algorithmic trajectories here.}
  \label{fig:conv_logistic_all}
\end{figure}

We experimented on the logistic regression problem \eqref{eq:logis} for three values for the dimension $d$: $d = 100, 400, 900.$ We computed the $C_d$ ratio statistics for a trial of the algorithm for these three dimensions, and display this information in Table \ref{tab:logistic_Cd}, which is based on the data from Figure \ref{fig:cd_plots_logistics}; the corresponding convergence behavior is shown in Figure \ref{fig:conv_logistic_all}. We observe that the maximum $C_d$ scales as around $2.5\sqrt{d}$, i.e., $\sO(\sqrt{d})$.

\begin{table}[ht]
  \captionsetup{width=\textwidth}
  \centering
  \caption{Statistics of the $C_d$ ratio with different dimensions $d$ for the logistic regression problem  \eqref{eq:logis}.}
  \label{tab:logistic_Cd}
  \begin{tabular*}{\textwidth}{@{\extracolsep{\fill}} c c c c @{}}
    \toprule[1.2pt]
    \textbf{Dimension $d$} & \textbf{Max} & \textbf{99th Percentile} & \textbf{Mean} \\
    \midrule
    $100$ & $6.9935$  & $2.8267$  & $0.5665$  \\
    $400$ & $47.2857$ & $24.1032$ & $3.7831$  \\
    $900$ & $75.4320$ & $64.0339$ & $16.4693$ \\
    \bottomrule[1.2pt]
  \end{tabular*}
\end{table}

\end{document}